\documentclass[reqno, 12pt]{amsart}

\usepackage{amssymb,amsmath}
\usepackage{amsthm}
\usepackage{upref}
\usepackage{enumerate}
\usepackage{amsfonts,bbold}

 \makeatletter
\def\LaTeX{\leavevmode L\raise.42ex
    \hbox{\kern-.3em\size{\sf@size}{0pt}\selectfont A}\kern-.15em\TeX}
 \makeatother
 
\DeclareMathOperator{\arccot}{arccot}

\numberwithin{equation}{section}

\newtheorem{lemma}{Lemma}[section]
\newtheorem{theorem}[lemma]{Theorem} 
\newtheorem{corollary}[lemma]{Corollary}
\newtheorem{proposition}[lemma]{Proposition}

\theoremstyle{definition}

\newtheorem{definition}[lemma]{Definition}

\newtheorem{remark}[lemma]{Remark}

\renewcommand{\det}{\operatorname{Det}}
\newcommand{\tr}{\operatorname{Tr}}

 \DeclareMathOperator*{\slim}{s-lim}

 \newcommand{\supp}{\operatorname{supp}}
  \newcommand{\e}{\eqref}

\newcommand{\q}{\quad}

\renewcommand{\d}{\delta}

\newcommand{\ov}{\overline}

    \newcommand{\Ch}{\operatorname{Ch}}

\renewcommand\Im{\operatorname{Im}}

\newenvironment{pf}{\begin{proof}}{\end{proof}}

\def\qqq{\mathrel{\subset\mkern-15mu\lower.38ex\hbox{${\scriptscriptstyle\rightarrow}$}}}

\let\cal\mathcal

\let\Bbb\mathbb

\begin{document}
\title [A point interaction for    the discrete Schr\"odinger operator]
 {A point interaction for    the discrete Schr\"odinger operator  and generalized Chebyshev polynomials}
\author{ D. R. Yafaev  }
\address{ IRMAR, Universit\'{e} de Rennes I, Campus de
  Beaulieu,  Rennes, 35042  FRANCE and SPGU, Univ. Nab. 7/9, Saint Petersburg, 199034 RUSSIA}
\email{yafaev@univ-rennes1.fr}
\subjclass[2000]{33C45, 39A70,  47A40, 47B39}

\date{17 May 2017}
 
  \keywords {Jacobi matrices, the discrete Schr\"odinger operator,  point interaction, resolvents,  explicit solutions, generalized Chebyshev polynomials}

\thanks{Supported by the grant RFBR 17-01-00668 A} 

\begin{abstract}
We consider  semi-infinite Jacobi matrices corresponding to a point interaction for  the discrete Schr\"odinger operator. Our goal is to find explicit expressions for the spectral measure,  the resolvent and other spectral characteristics of such Jacobi matrices. It turns out that the spectral analysis of this toy problem leads to a new class of orthogonal polynomials generalizing the classical  Chebyshev polynomials.
   \end{abstract}

\maketitle

\thispagestyle{empty}

\section{Introduction}

{\bf 1.1.}
As is well known, the theories of Jacobi operators given by three-diagonal matrices (see Subsection~2.1, for the precise definitions) and of differential operators $D a (x)D +b(x)$, $D=-i d/dx$,  are  to a large extent similar. This is true   for Jacobi operators acting in the space  $\ell^2 ({\Bbb Z})$ and   differential operators
acting in the space  $L^2 ({\Bbb R})$ as well as for the corresponding operators acting in the spaces  $\ell^2 ({\Bbb Z}_{+})$ and    $L^2 ({\Bbb R}_{+})$, respectively. We  refer to the book \cite{Teschl}  where this analogy is described in a very detailed way.
Both classes of the operators are very important in applications. For example, Jacobi operators play a substantial role in solid state physics (see, e.g., \S 1.5 of  \cite{Teschl}) while the Schr\"odinger operator is the basic object of quantum mechanics. Moreover, Jacobi operators   in the space  $\ell^2 ({\Bbb Z}_{+})$ are intimately related  (see, e.g.,  the classical book \cite{AKH}) to the theory of orthogonal polynomials. We refer to the books \cite{Ism, Sz} for all necessary information on  orthogonal polynomials.

We study  Jacobi operators    given in the space $\ell^2 ({\Bbb Z}_{+})$ by matrices
\begin{equation}
H_{a}=\frac{1}{2}
\begin{pmatrix}
 0&a& 0&0&0&\cdots \\
 a&0&1&0&0&\cdots \\
  0&1&0&1&0&\cdots \\
  0&0&1&0&1&\cdots \\
  \vdots&\vdots&\ddots&\ddots&\ddots&\cdots
\end{pmatrix} .
\label{eq:ZP+}\end{equation}
The parameter $a$ here is an arbitrary positive number. The operator $H_{1}$ plays the  role of the ``free"  differential operator $D^2$ in the space $L^2 ({\Bbb R}_{+})$ with the boundary condition $f(0)=0$. We consider $H_{a}$ as  a perturbation of the operator $H_{1}$  which is easy to analyze directly. Clearly,
  the perturbation $H_{a}-H_{1}$ has rank two  that allows us to find the resolvent of the operator $H_{a}$. The operator $H_{a}$ is a rare example where all spectral quantities, such as the spectral measure, eigenfunctions  (of the continuous spectrum), the wave operators, the scattering matrix, the spectral shift function, etc.,  can be calculated explicitly.

We show that
eigenfunctions of the operator $H_{a}$ are constructed in terms of a class of orthogonal polynomials ${\Ch}_{n}(z;a)$, $n=0,1,\ldots$,  generalizing the classical Chebyshev polynomials. It is very well-known that
  ${\Ch}_{n}(z;1)$ are the  Chebyshev polynomials
of the second kind.  It was also noted in the book \cite{Ber} that  ${\Ch}_{n}(z;\sqrt{2})$ are  the  Chebyshev polynomials of the first kind. We are studing  Jacobi matrices $H_{a}$ for  all $a>0$.  This naturally leads to a class of polynomials 
${\Ch}_{n}(z;a)$ parametrized by an arbitrary $a>0$.

The relation between Jacobi matrices and orthogonal polynomials is   the classical fact. In principle, it can be used 
in both ways. Let us mention the paper \cite{GS}  devoted to the inverse problem (of reconstruction of the Jacobi matrix by its spectral measure) where this is thoroughly discussed. The paper \cite{GS}  contains also numerous references to earlier works on this subject. 

We proceed from  the spectral  analysis of the Jacobi operator $H_{a}$ which can be performed in a very explicit way. 
  This leads to some formulas for the polynomials ${\Ch}_{n}(z; a)$. A part of these formulas is perhaps new even for the classical Chebyshev polynomials.

\medskip

{\bf 1.2.}
Compared to the continuous case, the operator $H_{a}$ plays the role of the self-adjoint realization $A_{\alpha}$ of
the differential operator $D^2$ in the space $L^2 ({\Bbb R}_{+})$ with a boundary condition $f'(0)=\alpha f(0)$ where $\alpha\in{\Bbb R}$ or of the  operator  $A^{(0)}$  with the boundary condition $  f(0)=0$. The operator $A^{(0)}$ is usually taken for the ``free'' operator, and the operators  $A_{\alpha}$  are interpreted as   Hamiltonians of a point interaction of two quantum particles. Of course the  essential spectra of the operators  $A_{\alpha}$ and $A^{(0)}$ coincide with  the half-axis $[0,\infty)$, and they coincide with the interval $[-1,1]$ for all operators $H_{a}$.
The study of the operators $H_{a}$ is more difficult and the formulas obtained are less trivial than for the operators $A_{\alpha}$. Actually, the operators $H_{a}$  are more close to self-adjoint realizations of the operator $D^2$ in the space ${\Bbb C}\oplus L^2 ({\Bbb R}_{+})$ (see the paper \cite{Y14x} or \S 4.7 of the book \cite{YA}).

Let us discuss the analogy between the operators 
  $A_{\alpha}$ and $H_{a}$ in  more details.     The role of $x\in  {\Bbb R}_{+}$ is played by the variable $n\in  {\Bbb Z}_{+}$, and the role of a function $f(x)$ is  played by a sequence
  $f=(f_{0},f_{1},\ldots)$. 
   The operator  $H_1$ acts by the formula $(H_{1}f)_{n}=2^{-1} (f_{n-1}+ f_{n+1})$ for all $n\in {\Bbb Z}_{+}$ if one imposes an artificial ``boundary condition" $f_{-1}=0$. So, the operator $H_{1}$ plays the role of $A^{(0)}$.  The operators $A_\alpha\geq 0$  for   $\alpha\geq 0$, and they have the    eigenvalues
   $-\alpha^{2}$ for all $\alpha<0$. The case $\alpha=0$ is critical: the operator $A_{0}\geq 0$, but its resolvent $(A_{0}-z I)^{-1}$ has a singularity at the bottom of its spectrum   (one says that the operator $A_{0}$ has a zero-energy resonance). The   operators $H_{a}$  have discrete eigenvalues (one below $-1$ and another one above $1$) if and only if $a>\sqrt{2}$.  The resolvent $(H_{\sqrt{2}}-z I)^{-1}$ of $H_{\sqrt{2}}$ has singularities at both edge points $\pm 1$ of its  spectrum. Thus  the operator $H_{\sqrt{2}}$  plays the role of $A_{0}$, and the values $a<\sqrt{2}$ (resp. $a>\sqrt{2}$) correspond to  the values $\alpha >0$ (resp. $\alpha <0$).
  
To a some extent, this paper was motivated by looking for a discrete analogue of the point interaction very well studied in the continuous case.  However, at a technical level, 
  the analogy between the discrete and continuous Schr\"odinger operators is not used  in this paper. Our approach relies on a direct calculation  of the resolvent of the operator $H_{a}$.
  
    In principle, our method applies to arbitrary finite rank perturbations of the operator $H_{1}$, but, for high ranks, the formulas become less explicit. In this context, we mention 
    the operator $  H_{a,b}$ which is obtained from  \e{eq:ZP+} if the first row of the matrix $  H_a$ is replaced by
    $(b,a,0,\cdots)$. The operator $  H_{a,b}$  is also a natural candidate for a description of the point interaction in the discrete case. Note that  self-adjoint realizations of the operator $D^2$ in the space ${\Bbb C}\oplus L^2 ({\Bbb R}_{+})$ are also parametrized by several numbers.
     The operator $  H_{a,b}-H_{1}$    has, in general,  rank three so that the operator $  H_{a,b}$ can   be easily studied. However, to keep formulas as simple as possible, we restrict our attention to  the operator $H_{a} =  H_{a,0}$. Note that the operator $H_{1,b}$ (and general diagonal rank one perturbations of $H_{1}$) was already well studied before (see, e.g., paper \cite{Gol}).  The methods of the orthogonal polynomials theory were used in this setting in
        \cite{BD}. 
 
\medskip

{\bf 1.3.}
The paper is organized as follows. Section~2 is of a preliminary nature. First, we briefly recall necessary facts about Jacobi matrices and their link to orthogonal polynomials. Then we study the operator $H_{1}$. The  results presented here are   well-known, but, to motivate our study of the operator $H_{a}$, we use somewhat different approach compared, for example,  to paper \cite{KS}. We explicitly relate the operator $H_{1}$ to an auxiliary operator ${\bf H}$ acting in the space $\ell^2({\Bbb Z} )$ by the formula $({\bf H}f)_{n} =2^{-1} (f_{n-1}+ f_{n+1})$. The operator ${\bf H}$  can be easily diagonalized by the discrete Fourier transform.

 Section~3 plays the central role. Here we find the resolvent $R_{a}(z)=(H_{a}-z I)^{-1}$ of the operator $H_{a}$ and its spectral family. We develop scattering theory for the pair $H_{1}$, $H_{a}$ in Section~4. In particular, we calculate the corresponding scattering matrix $S_{a}(\lambda)$ and the spectral shift function $\xi_{a}(\lambda)$. We also establish a link between eigenfunctions   of the operator $H_{a}$ and the wave operators for the pair $H_{1}$, $H_{a}$.  The traces $ \tr \big(H_{a} - H_{1} \big)$ and the moments of the spectral  measure of the operator $H_{a}$ are   calculated in Section~5.
  In  Section~6, we discuss a relation of our results on the operators $H_{a}$ with results on the corresponding Hankel operators. We also  study the limits $a\to 0$ and $a\to\infty$.  These limits turn out to be very  singular.  
  The operator $  H_{a,b}$  is briefly discussed in Section~7. Finally, we prove two simple, but general, results
  for arbitrary Jacobi operators in the Appendix.

  The author is grateful to L.~Golinskii for a useful correspondence and to the referee for the pertinent report.

 \section{Jacobi matrices. The discrete Schr\"odinger operator}

 {\bf 2.1.}
 Let us briefly recall some basic facts about Jacobi operators in   the space $\ell^2 ({\Bbb Z}_{+})$.
We denote by $e_{n}$, $ n\in {\Bbb Z}_{+}$,  the canonical  basis in this space, that is, all components of the vector $e_{n}$ are zeros, except the $n$-th component which equals $1$.
  Let $a_{0}, a_{1}, \ldots$, $b_{0}, b_{1}, \ldots$ be some real sequences and $a_{n} >0$ for all $n=0,1,\ldots$. Then the Jacobi operator $H$ is defined by the formula
 \begin{equation}
  H e_n= a_{n-1} e_{n-1} +b_{n} e_{n } + a_{n} e_{n+1}  , \q n\in {\Bbb Z}_{+} ,
\label{eq:J}\end{equation}
where we accept that $e_{-1}=0$. The operator $H$ is obviously symmetric. We assume that the sequences $\{a_{n}\}$ and $\{b_{n}\}$ are bounded, and hence     $H$ is a bounded operator in   the space $\ell^2 ({\Bbb Z}_{+})$.
Let us denote by $dE(\lambda)$  its spectral measure and put $d\rho(\lambda)=d(E(\lambda)e_{0}, e_{0})$.
Since the support of the measure $d\rho$ is bounded,  the set of all polynomials    is dense in the space $L^2 ({\Bbb R}; d\rho)$.

Let us define polynomials $P_{n}(z)$ by the recurrent relations
 \begin{equation}
z P_{n}(z)= a_{n-1} P_{n-1}(z) +b_{n} P_{n}(z) + a_{n} P_{n+1}(z)  , \q n\in {\Bbb Z}_{+} , \q z\in {\Bbb C}.
\label{eq:J1}\end{equation}
 We  accept that $P_{-1}(z)=0$  and put $P_0(z)=1$. Then $P_{n}(z)$  is a polynomial   of degree $n$ and its coefficient at $z^n$ equals $(a_{0}a_{1}\cdots a_{n-1})^{-1}$. Clearly,  the linear sets spanned by $\{1,z,\ldots,z^n\}$ and by $\{P_{0}(z),P_1(z), \ldots, P_n(z)\}$  coincide.

 Comparing relations \e{eq:J} and \e{eq:J1} and using recursion arguments, we see that
  \begin{equation}
e_{n} = P_{n}(H)e_{0}  
\label{eq:J2}\end{equation}
for all $n\in {\Bbb Z}_{+}$. Thus the set of all vectors $H^n e_{0}$, $n\in {\Bbb Z}_{+}$, is dense in $\ell^2 ({\Bbb Z}_{+})$, and hence   the spectrum of the operator $H$ is simple with $e_{0}$ being the generating vector.
It also follows from \e{eq:J2} that 
   \begin{equation}
d(E(\lambda)e_{n}, e_{m}) = P_{n}(\lambda) P_{m}(\lambda) d\rho(\lambda) 
\label{eq:J4}\end{equation}
whence
 \begin{equation}
\int_{-\infty}^\infty P_{n}(\lambda) P_{m}(\lambda) d\rho(\lambda) =\d_{n,m}
\label{eq:J5}\end{equation}
where $\d_{n,n}=1$ and $\d_{n,m}=0$ for $n\neq m$. Of course   $\{P_0 (\lambda),P_1 (\lambda),\ldots,P_{n} (\lambda)\}$ are obtained by the Gram-Schmidt orthogonalization of the set $\{1,\lambda,\ldots,\lambda^n\}$ in the space $L^2({\Bbb R}_{+}; d\rho)$.    

Let us now define a mapping $U: \ell^2 ({\Bbb Z}_{+})\to L^2 ({\Bbb R}; d\rho)$ by the formula 
 \begin{equation}
 (Ue_{n})(\lambda)=P_{n} (\lambda).
 \label{eq:Ua}\end{equation}
It is isometric according to \e{eq:J5}. It is also unitary because  the set of all polynomials  $P_n(\lambda) $, $n\in {\Bbb Z}_{+}$,  is dense in $L^2 ({\Bbb R}; d\rho)$. Finally, the intertwining property 
 \begin{equation}
(U H f) (\lambda)= \lambda (U f)(\lambda)
\label{eq:IP}\end{equation}
  holds. Indeed, it suffices to check it for $f=e_{n}$ when according to definition \e{eq:J}  $(U H e_{n}) (\lambda)$ coincides with the right-hand side of \e{eq:J} while $\lambda (U e_{n})(\lambda)$ equals its left-hand side.

We note also that \e{eq:J4} ensures the formula 
 \begin{equation}
\int_{-\infty}^\infty P_{n}(\lambda) P_{m}(\lambda) (\lambda-z)^{-1}d\rho(\lambda) =(R(z)e_{n}, e_{m}), \q \forall n,m\in {\Bbb Z}_{+},
\label{eq:J6}\end{equation}
which yields an expression for the integrals in the left-hand side provided the resolvent $R(z)= (H-zI)^{-1}$ is known.
Here and below $I$ is the identity operator.

Thus, for all sequence  $\{a_{n}\}$, $\{b_{n}\}$, one can construct the measure $d\rho(\lambda)$ such that the polynomials 
$P_{n} (\lambda)$ defined  by \e{eq:J1} are orthogonal in $L^2 ({\Bbb R}; d\rho)$. This assertion is known as the Favard theorem. Its standard proof presented here relies on the spectral theorem for self-adjoint operators and does not of course give an explicit expression for the measure $d\rho(\lambda)$. On the contrary, given a probability measure $d\rho(\lambda)$, one  can reconstruct $\{a_{n}\}$, $\{b_{n}\}$, that is,  the Jacobi matrix \e{eq:J} with the spectral measure $d\rho(\lambda)$. The solution of this (inverse) problem
is discussed from various points of view in the article \cite{GS}.

In this paper we study the case when $a_{0}=a/2$, $a_{n}=1/2$ for $n\geq 1$ and $b_{n}= 0$ for all $n\geq 0$.

 \medskip

  {\bf 2.2.}
  Let us first consider the  ``free" discrete Schr\"odinger operator (the infinite Jacobi matrix) ${\bf H}$  in the space     $\ell^2 ({\Bbb Z})$  given by the formula
\begin{equation}
{\bf H} e_{n}=\frac{1}{2} (e_{n-1}+e_{n+1})
\label{eq:Z}\end{equation}
where $e_{n}$, $n\in {\Bbb Z}$, is the canonical  basis in the space $\ell^2 ({\Bbb Z})$. We do not distinguish  in notation the scalar products in $\ell^2 ({\Bbb Z})$ and $\ell^2 ({\Bbb Z}_{+})$.
Evidently, the operator $\bf H$ can   be explicitly diagonalized. Indeed, let ${\cal F}$,
\begin{equation}
({\cal F} f) (\mu)= \sum_{n\in{\Bbb Z}}f_{n} \mu^n,  \q \mu \in {\Bbb T} , \q f= (\ldots,  f_{-1},f_{0}, f_{1},\dots) , \q f_{n}= (f,e_{n}), 
\label{eq:F}\end{equation}
 be the discrete Fourier transform. Let the unit circle
${\Bbb T} $ be endowed with the normalized Lebesgue measure 
\begin{equation}
d {\bf m}(\mu)= (2\pi i \mu)^{-1}d\mu, \q \mu\in {\Bbb T} .
\label{eq:id}\end{equation}
Then the operator ${\cal F}: \ell^2 ({\Bbb Z})\to L^2 ({\Bbb T})$  is unitary. Since
\begin{equation}
({\cal F} {\bf H} f) (\mu)=  \frac{\mu+\mu^{-1}}{2} ({\cal F} f) (\mu),
\label{eq:Z2}\end{equation}
  the spectrum of the operator ${\bf H} $ is absolutely continuous, has multiplicity $2$ and coincides with the interval $[-1,1]$.

Now it is    easy to calculate the resolvent  ${\bf R}(z)= ({\bf H}-z I)^{-1}$ of the operator ${\bf H}$.  Below we  fix the branch of the analytic function $\sqrt{z^2 -1}$ of $z\in {\Bbb C}\setminus [-1,1]$ by the condition $\sqrt{z^2 -1}>0$ for $z>1$. Then it equals $\pm i\sqrt{1-\lambda^2}$ for $z=\lambda\pm i0$, $\lambda\in (-1,1)$, and $\sqrt{z^2 -1}< 0$ for $z< -1$. 

 \begin{lemma}\label{RZ}
  For all $n,m\in{\Bbb Z}$, we have
 \begin{equation}
({\bf R}(z)e_{n}, e_{m})=  -\frac{(z-\sqrt{z^2 -1})^{|n-m|} }{\sqrt{z^2 -1}}.
\label{eq:Z4}\end{equation}
 \end{lemma}

 \begin{pf}
  It follows from  \e{eq:F} -- \e{eq:Z2} that
\begin{align}
({\bf R}(z)e_{n}, e_{m})&= (2\pi i)^{-1}\int_{\Bbb T} \big(\frac{\mu+\mu^{-1}}{2}-z\big)^{-1}\mu^{n-m-1} d\mu
\nonumber\\
&=(\pi i)^{-1}\int_{\Bbb T} \big( \mu^2 -2z\mu +1\big)^{-1}\mu^{n-m} d\mu,\q z\in {\Bbb C}\setminus [-1,1].
\label{eq:Z3}\end{align}
The  equation $ \mu^2 -2z\mu +1=0$  has the roots 
 \begin{equation}
 \mu_{\pm}=z\pm \sqrt{z^2 -1}.
 \label{eq:Z3a}\end{equation}
  By  the proof of \e{eq:Z4}, we  may suppose that $z>1$.
  Then $\mu_{+} >1$  and $\mu_{+} \mu_- =1$. Since ${\bf H} $ commutes with the complex conjugation, we have
  \[
  ({\bf R}(z)e_{n}, e_{m})=  ({\bf R}(z)e_m, e_n),
  \]
so that  we can suppose $n\geq m$ in \e{eq:Z3}. Then $\mu_{-}$ in the only singular point of the integrand  in \e{eq:Z3} inside the unit circle. Calculating the residue at this point, we find that
\[
({\bf R}(z)e_{n}, e_{m})=  2 (\mu_{-}-\mu_{+})^{-1}  \mu_{-}^{n-m} .
\]
Substituting here expressions \e{eq:Z3a}, we arrive at \e{eq:Z4}.
 \end{pf}

{\bf 2.3.}
Now we are in  a position to study the discrete  ``free" Schr\"odinger operator  (the semi-infinite Jacobi matrix) $ H_{+}$  in   the space $\ell^2 ({\Bbb Z}_{+})$ with elements  $f= ( f_{0}, f_{1},\dots)  $.  It  is defined by the formula
\begin{equation}
H_{+} e_0=\frac{1}{2} e_{ 1} , \q  H_{+} e_{n}=\frac{1}{2} (e_{n-1}+e_{n+1}), \q n\geq 1.
\label{eq:Z+}\end{equation}
Of course the operator $ H_{+}$  is given by  matrix \e{eq:ZP+} where $a=1$.

Our first goal is to calculate the resolvent $R_{+}(z)= (H_{+} -zI )^{-1}$ of the operator   $ H_{+}$. To that end, we introduce an auxiliary operator $H_{-}$ in the space $\ell^2 ({\Bbb Z}_{-})$ where ${\Bbb Z}_{-}={\Bbb Z} \setminus {\Bbb Z}_{+}$
by the formula
\begin{equation}
H_{-} e_{-1}=\frac{1}{2} e_{ -2} , \q  H_{-} e_{-n}=\frac{1}{2} (e_{-n-1}+e_{-n+1}), \q n\geq 2.
\label{eq:Z-}\end{equation}
The operators $H_{+}$ and $H_{-}$ are of course unitarily equivalent, but we do not need this fact. 
 Let us   consider
$H=H_{+}\oplus H_{-}$ as an operator in   the space $\ell^2 ({\Bbb Z} )=\ell^2 ({\Bbb Z}_{+})\oplus \ell^2 ({\Bbb Z}_-)$. Comparing formulas \e{eq:Z}, \e{eq:Z+} and \e{eq:Z-}, we see that the difference
\begin{equation}
V=H-{\bf H}= -\frac{1}{2}  (\cdot, e_{0}) e_{-1}  -\frac{1}{2}  (\cdot, e_{-1}) e_0
\label{eq:V}\end{equation}
is an operator in $\ell^2 ({\Bbb Z})$ of rank $2$. This allows us to easily obtain an explicit expression for the resolvent $R(z) =(H-z I)^{-1}$ of the operator $H$.

Let us set 
\[
T(z)=V-V R(z) V.
\]
In terms of the operator $T(z)$, the resolvent equation for the pair ${\bf H}$, $H$ can be written as
\begin{equation}
T(z)=V-V {\bf R}(z) T(z).
\label{eq:T1}\end{equation}

Denote
\begin{equation}
 \omega (z) =z-\sqrt{z^2 -1}= (z + \sqrt{z^2 -1})^{-1}.
\label{eq:ome}\end{equation}
Below we often use the obvious identities:
\begin{equation}
 1-\omega(z)^2
=2\sqrt{z^2-1}\,\omega(z), \q  1+\omega(z)^2
=2z \omega(z).
\label{eq:id2}\end{equation}

\begin{lemma}\label{TZ}
The solution $T(z) $ of the equation  \e{eq:T1} is given by the formula
\begin{equation}
2 T(z)f=   (-f_{0} +  f_{-1}\omega (z) ) e_{-1} +  (-f_{-1} +  f_{0}\omega (z) )  e_0,\q f_{n} =(f,e_{n}) .
\label{eq:T6}\end{equation}
 \end{lemma}
 
\begin{pf}
In view of expression \e{eq:V},  equation \e{eq:T1} implies that
\begin{equation}
2 T(z)f= -  f_{0}  e_{-1}  - f_{-1} e_0
+    (  {\bf R}(z) T(z)f, e_{0}) e_{-1}  +   ( {\bf R}(z) T(z)f,   e_{-1}) e_0 
\label{eq:T2}\end{equation}
where $ f_{n} =(f, e_n)$
for all $ f\in   \ell^2 ({\Bbb Z})$. Let us take the scalar products of this equation with $ {\bf R}(\bar{z})e_{-1}$ and $ {\bf R}(\bar{z})e_0$. This leads to a system of two equations for $u_{-1} (z)=({\bf R}(z)T(z)f,  e_{-1})$ and $u_0 (z)=({\bf R}(z)T(z)f,  e_{0})$:
\begin{equation}
\left.
\begin{aligned}
2u_{-1} (z) = &-   ({\bf R}(z) e_{-1} ,e_{-1} ) f_{0} -  ({\bf R}(z) e_0 ,e_{-1} ) f_{-1}
\\ & +     ({\bf R}(z) e_{-1} ,e_{-1} )   u_{0} (z ) +   ({\bf R}(z) e_{0} ,e_{-1} ) u_{-1} (z) 
\\
2u_{0} (z) = &-     ({\bf R}(z) e_{-1} ,e_{0} ) f_{0} -   ({\bf R}(z) e_0 ,e_{0} ) f_{-1}
\\
&+       ({\bf R}(z) e_{-1} ,e_{0} ) u_{0} (z) +    ({\bf R}(z) e_{0} ,e_{0} ) u_{-1} (z)
\end{aligned}
\right\} .
\label{eq:T3}\end{equation}

To solve   system  \e{eq:T3}, we   take into account expression
\e{eq:Z4}. Then \e{eq:T3}  can be rewritten as
\[
\left.
\begin{aligned}
(2\sqrt{z^2 -1} +\omega (z) ) u_{-1}  (z)  + u_{0}  (z) = & f_{0}  + f_{-1}\omega  (z) 
\\
u_{-1}  (z)  + (2\sqrt{z^2 -1} +\omega  (z)  )     u_{0}  (z) = & f_{0} \omega  (z)  + f_{-1} 
\end{aligned}
\right\}
\]
which yields
\[
 u_{-1} (z)=   f_{0}  \omega (z), \q  u_{0} (z)=   f_{-1}\omega (z).
\]
Therefore the representation \e{eq:T6} is  a direct consequence of equation \e{eq:T2}.
\end{pf}

Now we are in a position to calculate matrix elements of the resolvent
\begin{equation}
R (z)={\bf R}(z)- {\bf R}(z) T(z) {\bf R}(z).
\label{eq:T7}\end{equation}
 
 \begin{proposition}\label{RZ+}
Let the operator $H_{+}$ be defined in   the space $\ell^2 ({\Bbb Z}_{+})$  by relations
\e{eq:Z+}. Then the matrix elements of its resolvent $R_{+} (z) =(H_{+} -z I)^{-1}$ are given by the formula
 \begin{equation}
(R_{+}(z)e_{n}, e_{m})=  \frac{(z-\sqrt{z^2 -1})^{n+m+2}-(z-\sqrt{z^2 -1})^{|n-m|} }{\sqrt{z^2 -1}}=:r_{n,m}(z)
\label{eq:Z4+}\end{equation}
 for all $n,m\in{\Bbb Z}_{+}$. In  particular,
 \begin{equation}
(R_{+}(z)e_0, e_0)= 2 (  \sqrt{z^2 -1}-z) .
\label{eq:Z40}\end{equation}
 \end{proposition}
 
  \begin{pf}
   It follows from formula \e{eq:Z4} and Lemma~\ref{TZ} that, for all $n,m\geq 0$,  
   \begin{equation}
(  {\bf R}(z)  T(z)  {\bf R}(z) e_{n} ,e_{m}) = \frac{\omega(z)^{n+m}}{z^2-1} (T(z) (\omega(z)e_{-1}+e_{0}) ,\omega(\bar{z})e_{-1}+e_{0}) .
\label{eq:TS}\end{equation}
According to \e{eq:T6} and the first identity \e{eq:id2}, we have
\[
2(T(z) (\omega(z)e_{-1}+e_{0}) ,\omega(\bar{z})e_{-1}+e_{0})=\omega(z)^3-\omega(z)=-2\sqrt{z^2-1}\omega(z)^2
\]
and therefore \e{eq:TS} implies
\begin{equation}
(  {\bf R}(z)  T(z)  {\bf R}(z) e_{n} ,e_{m}) = - \frac{\omega(z)^{n+m+2}}{\sqrt{z^2-1} }   .
\label{eq:TS1}\end{equation}
Let us now take into account representation \e{eq:T7} and observe that
$(R_{+} (z)e_{n}, e_{m}) = (R (z)e_{n}, e_{m})$ if $n\geq 0$ and $m\geq 0$. Thus  
 formula  \e{eq:Z4+} is an immediate corollary of \e{eq:Z4} and \e{eq:TS1}.
   \end{pf}
 
The following result is a direct consequence  of    the standard relation 
  \[
2\pi i \frac{d(E_{+}(\lambda)e_0, e_0)} {d\lambda}=  (R_{+}(\lambda+ i0)e_0, e_0)-(R_{+}(\lambda- i0)e_0, e_0)
\]
 between the spectral measure  $d E_{+}(\lambda)$ of the self-adjoint operator $H_{+}$ and its resolvent.

 \begin{corollary}\label{EZ+}
 For all $\lambda\in (-1,1)$, we have
 \begin{equation}
d(E_{+}(\lambda)e_0, e_0)= 2 \pi^{-1}  \sqrt{1-\lambda^2 }d\lambda. 
\label{eq:ZE+}\end{equation}
 \end{corollary}

\section{Point interaction. The spectral measure}

  {\bf 3.1.}
 Here we consider a   generalization $H_{a}$ of the operator \e{eq:Z+}  given in   the space $\ell^2 ({\Bbb Z}_{+})$ by the matrix
 \e{eq:ZP+}. To put it differently, this operator is defined by the formulas
\begin{equation}
H_a e_0=\frac{a}{2} e_{ 1} , \q  H_a e_1=\frac{1}{2} (a e_0 +e_2), \q  H_a e_{n}=\frac{1}{2} (e_{n-1}+e_{n+1}), \q n\geq 2,
\label{eq:ZP}\end{equation}
where $a>0$.
Obviously,   $  H_{1}=H_{+}$. Of course, the operators $H_{a}$ are particular cases of the Jacobi operators discussed in Subsection~2.1 corresponding to the case $a_{0}=a/2$, $a_{n}=1/2$ for $n\geq 1$ and $b_{n}=0$ for all $n\geq 0$. So, all the results exposed there are automatically true now.

Our first goal here is to give explicit expressions for  the objects introduced in Subsection~2.1. For the operator $H_{a}$, they will be denoted $d E_{a} (\lambda)$, $d \rho_{a}(\lambda)$,    $R_{a} (z) =(H_{a}-zI)^{-1}$, $U_{a}$, etc.
 Let us start with the polynomials $P_{n} (z)$. For the case of  the operator $H_{a}$, they will be denoted $\Ch_{n} (z;a)$.
   To be precise, we  accept the following

  \begin{definition}\label{CheR}
 The polynomials   $\Ch_{n} (z;a)$ are defined
 by the recurrent relations: 
 $\Ch_{0} (z ;a)=1$,   $\Ch_1 (z;a)=2a^{-1} z$, 
 \begin{equation}
a \Ch_0(z;a)+ \Ch_2 (z;a) =2 z \Ch_1  (z;a) 
\label{eq:Chab}\end{equation}
and
\begin{equation}
 \Ch_{n-1} (z;a)+ \Ch_{n+1} (z;a) =2 z \Ch_{n}  (z;a), \q n\geq 2 .
\label{eq:chab}\end{equation}
\end{definition}

Note that $\Ch_{n}  (z;1)={\sf U}_{n}(z)$ (the Chebyshev polynomials   of the second kind)  and $\Ch_{n}  (z;\sqrt{2})={\sf T}_{n}(z)$  (the Chebyshev polynomials   of the first kind) for all $n\in {\Bbb Z}_{+}$. For an arbitrary $a>0$, we use the term ``generalized  Chebyshev polynomials" for $\Ch_{n}  (z;a)$.  It is possible to give an explicit expression for these polynomials.

 \begin{proposition}\label{Cheb}
Let us set
\begin{equation}
\gamma_{\pm}  (z;a)= \frac{a}{2} \pm  \frac{a^2-2}{2a}\frac{z}{\sqrt{z^2-1}}.
\label{eq:CH1}\end{equation}
Then
\begin{equation}
  \Ch_{n}  (z;a) =\gamma_{+} (z;a)\omega(z)^n + \gamma_{-} (z;a)\omega(z)^{-n},\q \forall n\geq 1,
\label{eq:CH2}\end{equation}
where the function $\omega(z)$ is defined by formula \e{eq:ome}.
\end{proposition}
 
 \begin{pf}
 Since
 \[
 \omega(z)+ \omega(z)^{-1} =2z,
 \]
 both sequences $f_{n}^{(+)}(z)= \omega(z)^n$ and $f_{n}^{(-)}(z)= \omega(z)^{-n}$ satisfy the equations 
 \[
 f_{n-1}(z)+ f_{n+1}(z)=2z f_{n}(z),\q n\geq 2.
 \]
 Therefore their arbitrary linear combination \e{eq:CH2}   solves   equations \e{eq:chab}.  To find the constants $\gamma_{+} (z;a)$ and $\gamma_{-} (z;a)$, it  remains to satisfy  the equations $  \Ch_1  (z;a) = 2a^{-1}  z$ and \e{eq:Chab} where $\Ch_0  (z;a)=1$. This yields the system
 \[
\left.
\begin{aligned}
 \gamma_{+} (z;a) \omega(z)+ \gamma_{-} (z;a) \omega(z)^{-1}&= 2a^{-1}  z
\\
 a+\big( \gamma_{+} (z;a) \omega(z)^2 + \gamma_{-} (z;a) \omega(z)^{-2}\big)&= 4a^{-1}  z^2
\end{aligned}
\right\}.
\]
  It is easy to see that its solution   is given by formula \e{eq:CH1}.
  \end{pf}
  
   \begin{remark}\label{Cher}
 All  functions $\gamma_\pm (z;a)$ and $\omega(z)$ are holomorphic on the set ${\Bbb C}\setminus [-1,1]$ only.
   Nevertheless their combination in the right-hand side of \e{eq:CH2} is a polynomial.
        \end{remark}
  
   \begin{remark}\label{Chepm}
    Representation \e{eq:CH2} remains true for $z=\pm 1$. In this case 
 \begin{equation}
    \Ch_{n}  (\pm 1;a) =(\pm 1)^n \big(   2n a^{-1} -(n-1)a \big),\q n\geq 1.
 \label{eq:CH3}\end{equation}
    Indeed, let $z\to\pm 1$ in \e{eq:CH2}. Observe that  $\omega(\pm 1)=\pm 1$,
    \[
    \omega(z)^{n} - \omega(z)^{-n}= - 2n z^{n-1 }\sqrt{z^2-1}+  O(z^2-1)
    \]
    and
    \[
    \gamma_{+} (z)= \frac{a^2 - 2}{2a}\frac{z}{\sqrt{z^2-1}}+O(1)\q {\rm as}  \q z\to\pm 1.
    \]
       Since  $ \gamma_{+} (z;a)  + \gamma_{-} (z;a)=a$, formula \e{eq:CH3} follows from \e{eq:CH2} in the limit $z\to\pm 1$.
       \end{remark}

  Let us now use \e{eq:CH2}  for $z=\lambda+ i0$, $ \lambda\in (-1,1)$, and observe that according to
 \e{eq:ome} and \e{eq:CH1} 
   \begin{equation}
\omega (\lambda \pm i0)= \lambda \mp i{\sqrt{1-\lambda^2}}
\label{eq:ome1}\end{equation}
and
\begin{equation}
\gamma_{\pm}  (\lambda + i0;a)= \frac{a}{2} \mp i  \frac{a^2-2}{2a}\frac{\lambda}{\sqrt{1-\lambda^2}}.
\label{eq:CH12}\end{equation}
Note that $  \Ch_{n}  (-\lambda;a)= (-1)^n \Ch_{n}  (\lambda;a)$ because $\omega (-\lambda + i0) =- \ov{\omega (\lambda + i0)}$ and $\gamma_{+}  (-\lambda + i0;a) = \gamma_{-}  (\lambda + i0;a)$.

  Proposition~~\ref{Cheb}  implies the following assertion.  
  
\begin{corollary}\label{Cheb11}
Let $\lambda\in (-1,1)$. Put $\lambda= \cos \theta$ where $\theta=\theta (\lambda)\in (0,\pi)$ and
\begin{equation}
 \varkappa (\theta;a)= 
 \frac{\sqrt{a^4 +4 (1-a^2)\cos^2\theta}}{a\sin \theta},\q   \d(\theta;a)= \arctan \Big(\frac{2-a^2}{a^2} \cot \theta\Big).
\label{eq:CH13}\end{equation}
Then  
\begin{equation}
  \Ch_{n}  (\lambda;a) = \varkappa (\theta;a)  \cos \big( n\theta- \d(\theta;a)\big) ,\q   n\geq 1 .
\label{eq:CHx}\end{equation}
 \end{corollary}
 
 \begin{pf}
With definitions \e{eq:CH13}, we have
 \[
 \omega (\lambda +  i0)=e^{-i   \theta }, \q   \gamma_{+}  (\lambda + i0;a)= 2^{-1}\varkappa (\theta;a) 
 e^{ i  \d(\theta;a)}.
 \]
 Therefore \e{eq:CHx} is a direct consequence of \e{eq:CH2}  for $z=\lambda+ i0$.
  \end{pf}

  Representation \e{eq:CHx} extends the well-known formula (10.11.2) in the book \cite{BE}   for the classical Chebyshev polynomials to  arbitrary $a>0$.

       Of course passing in \e{eq:CHx} to the limit $\lambda\to\pm 1$, we recover relation      \e{eq:CH3}.

\bigskip

 {\bf 3.2.}
 Now we perform an explicit spectral analysis of the operators $H_{a}$ for all $a>0$. In particular, we find the orthogonality measure $d\rho_{a} (\lambda)= d(E_{a} (\lambda)e_{0}, e_{0})$ for the polynomials $\Ch_{n}  (\lambda;a)$. 
 Let us consider the operator $H_{a}$ as a perturbation of $H_{1}$. Obviously,  the perturbation
\begin{equation}
V_{a}= H_a-H_{1}= \alpha(\cdot, e_{0}) e_{1}  + \alpha (\cdot, e_{1}) e_0, \q\alpha= \frac{a-1}{2} ,
\label{eq:Va}\end{equation}
has rank $2$. Of course, the essential spectrum of $H_{a}$  is the same as that of $H_{1}$, i.e., it coincides with the interval $[-1,1]$. 

We need a particular case of Proposition~\ref{RZ+}.  It can be easily deduced from formula \e{eq:Z4+}  if the identities \e{eq:id2}  are taken into account.

 \begin{lemma}\label{RZr}
 Let the function $\omega(z)$ be  defined by relation \e{eq:ome}, and let
 \begin{equation}
r_{n,m}(z)=(R_1 (z)e_{n}, e_{m}) .
\label{eq:Vr}\end{equation}
 Then
\begin{equation}
\begin{split}
r_{n,0}(z)&=-2\omega^{n+1} (z), \q n\geq 0,
\\
r_{n,1}(z)& =  -4z \omega^{n+1}(z),\q n\geq 1.
 \end{split}
\label{eq:omr}\end{equation}
 \end{lemma}
 
 Our first goal is to calculate the perturbation determinant
\begin{equation}
D_{a}(z)=\det \big(I+ V_{a} R_{1}(z)\big)
\label{eq:PD}\end{equation}
for the pair $H_{1}$, $H_{a}$.

 \begin{proposition}\label{PD}
For all $a>0$, the perturbation determinant is given by the formula
\begin{equation}
D_{a}(z)=1+  ( 1-a^2) \omega (z)^2.  
\label{eq:D}\end{equation}
In particular,
\[
D_{\sqrt{2}}(z)=2 \sqrt{z^{2}-1} \, \omega (z).
\]
 \end{proposition}
 
 \begin{pf}
 It follows from definitions  \e{eq:Va} and  \e{eq:Vr} that
\begin{multline*}
D_{a}(z)=\det  \begin{pmatrix}
1+\alpha r_{1,0} (z)&\alpha r_{0,0} (z) \\
\alpha r_{1,1} (z)&1+\alpha r_{0,1} (z) 
\end{pmatrix}
\\
=(1+\alpha r_{1,0} (z)) (1+\alpha r_{0,1} (z)) -\alpha^2  r_{0,0} (z)  r_{1,1} (z) .
\end{multline*}
Substituting here expressions \e{eq:omr}, we get \e{eq:D}.
 \end{pf}
 
 Note that formulas more general than \e{eq:D} were obtained earlier in \cite{Gol}.
 
 Next, we calculate eigenvalues of the operator $H_{a}$.
 
  \begin{proposition}\label{RZE}
If $a\leq\sqrt{2}$, then the set $(-\infty, -1)\cup (1,\infty)$ consists of regular points of  the operator $H_{a}$. 
If $a>\sqrt{2}$, then    the operator $H_{a}$ has exactly two isolated eigenvalues
\begin{equation}
\lambda_{\pm} (a)=\pm \frac{a^2}{2\sqrt{a^2-1}  }.
\label{eq:eig}\end{equation}
The points $1$ and $-1$ are not eigenvalues of the operators $H_{a}$.
 \end{proposition}
 
   \begin{pf}
   Suppose first that $ |\lambda|>1$. Recall that $\lambda$ is an eigenvalue of   $H_{a}$ if and only if $D_{a}(\lambda)=0$.
   In view of relation \e{eq:D} and the second identity \e{eq:id2} this equation is equivalent to $2\lambda=a^2 \omega(\lambda)$ or, by definition \e{eq:ome},
   \[
(a^2-2) \lambda= \pm a^2\sqrt{\lambda^2 -1} \q{\rm if}\q \pm \lambda>1.
\]
     Obviously, this equation has solutions if and only if $a^2>2$. In this case these  solutions are given by formula \e{eq:eig}.
     
     Let now $\lambda=  1$ or $\lambda=  -1$ and $(H_{1}+V_{a})\psi=\lambda\psi$ for some $\psi\in \ell^2 ({\Bbb Z}_{+})$. Put  $\psi_{n}= (\psi, e_n)$ and
       \[
\varphi= (H_{1}-\lambda)\psi =- V_{a}\psi .
\]
It follows from formula \e{eq:Va} that
 \[
\varphi =-\alpha \psi_{0}e_1 -\alpha \psi_1 e_{0}
\]
and hence
\[
\psi= R_{1} (\lambda)\varphi =-\alpha \psi_{0}R_{1} (\lambda)e_1 -\alpha \psi_1 R_{1} (\lambda)e_{0}.
\]
Using formulas \e{eq:omr} where $ \omega(\lambda)=\lambda$, we find that
\[
\psi_{n}=  2\alpha   \lambda^{n+1}(2\lambda \psi_{0}+\psi_{1}), \q n\geq 0.
\]
Comparing these two equations for $n=0$ and $n=1$, we see that $\psi_{1}=\lambda \psi_{0}$ whence
$\psi_{n}= 6\alpha \lambda^n \psi_{0}$. Since $ \lambda=\pm 1$, the sequence
  $\psi\in \ell^2 ({\Bbb Z}_{+})$ for $\psi=0$ only.
 \end{pf}

\bigskip

 {\bf 3.3.}
Now we are in a position  to calculate the resolvent $R_{a} (z)=(H_{a} -z I)^{-1}$ of the operator $H_{a}$.   Similarly to Subsection~2.3, it is more convenient to work with the operator
\begin{equation}
T_{a}(z)=V_{a}-V_{a}R_a (z) V_{a}.
\label{eq:Ta}\end{equation}
Then
\begin{equation}
R_{a} (z)=R_{1}(z)- R_{1}(z) T_{a}(z) R_{1}(z)
\label{eq:T7A}\end{equation}
and the resolvent equation for the pair $H_1$, $H_{a}$ can be written as
\begin{equation}
T_{a}(z)=V_{a}-V_{a} R_{1}(z) T_{a}(z).
\label{eq:Ta1}\end{equation}

 \begin{proposition}\label{TZa}
Let the operator $H_{a}$ be defined by formulas \e{eq:ZP}. Then for all $f \in \ell^2 ({\Bbb Z}_{+})$, we have
\begin{equation}
2 D_{a}(z) T_{a}(z)f=   (t_{0,0} (z) f_{0}+ t_{0,1} (z) f_{1})e_{0}+   (t_{1,0} (z) f_{0}+ t_{1,1} (z) f_{1})e_1
\label{eq:Ta6}\end{equation}
where   the perturbation determinant $D_{a}(z) $ is given by formula  \e{eq:D},  $f_{n}=(f, e_{n})$ and
\begin{equation}
\begin{split}
t_{0,0} (z) =  2 (a-1)^2  z \omega(z)^2,\q & t_{1,1} (z) =(a-1)^2 \omega(z),
\\
t_{0,1} (z) = t_{1,0} (z) = a- 1- & (a-1)^2 \omega(z)^2 .
\end{split}
\label{eq:Ta7}\end{equation}
 \end{proposition}

 \begin{pf}
In view of expression \e{eq:Va}  for $V_{a}$,    equation \e{eq:Ta1} can be rewritten as
\begin{equation}
T_{a}(z)f=  \alpha  f_{1} e_0 +\alpha  f_{0} e_{1}  -\alpha  (R_{1}(z) T_{a}(z)f,  e_1) e_0 
-\alpha  ( R_{1}(z) T_{a}(z)f,  e_{0}) e_{1}  .
\label{eq:Ta2}\end{equation}
 Let us take the scalar products of this equation with $R_{1}(\bar{z})e_0$ and $R_{1}(\bar{z})e_1$. This leads to a system of two equations for  $u_0 (z)=(R_{1}(z) T_{a}(z)f,  e_{0})$ and $u_{1} (z)=(R_{1}(z) T_{a}(z)f,  e_{1})$:
\[
\left.
\begin{aligned}
u_0 (z) = & \alpha     (R_{1}(z) e_1,e_0 )f_{0} + \alpha   (R_{1}(z) e_0 ,e_0 ) f_{1}
\\ & -\alpha     (R_{1}(z) e_1 ,e_0 ) u_{0} (z)  -\alpha    (R_{1}(z) e_{0} ,e_0 ) u_{1} (z)
\\
u_1 (z) = & \alpha    (R_{1}(z) e_1,e_1 )  f_{0}+ \alpha   (R_{1}(z) e_0 ,e_1 ) f_{1}
\\ & -\alpha      (R_{1}(z) e_1 ,e_1 ) u_{0} (z) -\alpha    (R_{1}(z) e_{0} ,e_1 ) u_{1} (z)
\end{aligned}
\right\}.
\]

Let us solve this system. Using notation \e{eq:Z4+}, we can rewrite it as
\begin{equation}
\left.
\begin{aligned}
 (1+\alpha r_{1,0}(z))u_{0} (z) +\alpha r_{0,0} (z) u_1 (z)= & \alpha r_{1,0} (z) f_{0}  + \alpha r_{0,0} (z) f_{1} 
\\
 \alpha r_{1,1}(z)u_{0} (z) + (1+\alpha r_{0,1}(z)) u_1(z)= & \alpha r_{1,1} (z) f_{0}  + \alpha r_{0,1} (z) f_{1} 
\end{aligned}
\right\}.
\label{eq:Ta4}\end{equation}
 It can be easily checked that the solution of system \e{eq:Ta4} is given by the formulas
\begin{equation}
\begin{aligned}
D_{a} (z)  u_0 (z) = & -a (a-1)\omega(z)^2  f_{0}  - (a-1) \omega(z) f_{1} ,
\\
D_{a} (z) u_1 (z)= &  -2 (a-1) z \omega(z)^2 f_0  - a (a-1)\omega(z)^2  f_{1} .
\end{aligned}
\label{eq:Ta5}\end{equation}
According to \e{eq:Ta2}
 we have 
 \[
T_{a}(z)f= \alpha ( f_1 -u_1 (z)) e_{0} +\alpha ( f_0 -u_{0} (z)) e_{1},
\]
which in view of \e{eq:Ta5} yields the representation \e{eq:Ta6}, \e{eq:Ta7}.
 \end{pf}

To find the spectral measure of the operator $H_{a}$, we 
  calculate the matrix element $ (R_{a} (z) e_{0} , e_{0})$ of the resolvent; it is also known as the Weyl $m$-function. Other matrix elements of $R_{a} (z)$ will be found in Subsection~5.1.
  
  \begin{theorem}\label{TTab}
  The representation
\begin{equation}
( R_{a} (z)  e_0,e_{0})=-\frac{2\omega(z)}{ D_{a}(z)} =\frac{2 \omega(z)}{( a^2-1)\omega(z)^2 -1}=\frac{2 }{( a^2-2)z- a^2\sqrt{z^2-1}}
\label{eq:Ta7X}\end{equation}
holds.
 \end{theorem}

  \begin{pf}
   In view of the identity \e{eq:id2},
 it follows from formula \e{eq:Z4+} that $(R_{1}(z) e_0, e_{0})=-2 \omega(z)$ and $(R_{1}(z) e_0, e_1)=-2 \omega(z)^2$.
 Therefore representation  \e{eq:Ta6} implies that
\begin{align*}
2D_{a}(z)(R_{1}(z)  T_{a} (z)R_{1}(z) e_0, e_{0}) =& 8 D_{a}(z)\omega(z)^2  (  T_{a} (z) (e_0+\omega(z) e_{1}), e_0+\omega(\bar{z}) e_{1})
\nonumber\\
= & 4 \omega(z)^2  ( t_{0,0}(z)+2\omega(z)t_{0,1}(z)+\omega(z)^2 t_{1,1}(z)).
\end{align*}
Substituting here expressions \e{eq:Ta7}, we see that
\[
2D_{a}(z)(R_{1}(z)  T_{a} (z)R_{1}(z) e_0, e_{0}) = 4 (a^2-1) \omega(z)^3.
\]
 In view of \e{eq:T7A}, this yields  the representation
\[
( R_{a} (z)  e_0,e_{0})-( R_1 (z)  e_0,e_{0})=- 2 (a^2-1)\frac{\omega(z)^3 }{ D_{a}(z)} .
\]
Using also equalities \e{eq:Z40} for $( R_1 (z)  e_0,e_{0})$ and  \e{eq:D} for $ D_{a}(z)$, we get the first equality \e{eq:Ta7X}. To get other equalities, we use again  \e{eq:D} and  \e{eq:ome}.
     \end{pf}

 Observe now that the denominator in \e{eq:Ta7X} is a continuous function of $z\in{\Bbb C}\setminus [-1,1]$ (up to the cut along $[-1,1]$),  and it is not $0$ for $z=\lambda\pm i0$, $\lambda\in (-1,1)$. Applying the formula  
  \[
2\pi i \frac{d(E_{a}(\lambda)e_0, e_0)} {d\lambda}=  (R_{a}(\lambda+ i0)e_0, e_0)-  (R_{a}(\lambda - i0)e_0, e_0),
\]
we find the spectral measure $d (E_{a}(\lambda) e_{0}, e_{0})$ of the self-adjoint operator $H_{a}$ for $\lambda\in (-1,1)$. Finally, calculating the residues of function \e{eq:Ta7X} at the points $\lambda_{\pm} (a)$, we obtain the spectral measure at these points.

 \begin{theorem}\label{TZab}
Let the Jacobi operator $H_{a}$ in the space $\ell^2 ({\Bbb Z}_{+})$ be defined by formula \e{eq:ZP+}.
 Then:
  \begin{enumerate}[\rm(i)]

\item
 The spectrum of the operator $H_{a}$ is absolutely continuous on $[-1,1]$ and
 \begin{equation}
d(E_{a}(\lambda)e_0, e_0)= 2 \pi^{-1} \frac{a^2\sqrt{1-\lambda^2 }}{a^4-4 (a^2-1) \lambda^2}d\lambda =: d\rho_{a}(\lambda) 
\label{eq:ZEa}\end{equation}
for all $\lambda\in (-1,1)$.

\item
If $a\leq\sqrt{2}$, then the set $(-\infty,-1)\cap (1,\infty)$ consists of regular points of the operator $H_{a}$.
If $a>\sqrt{2}$, then the operator $H_{a}$ has two simple eigenvalues  $ \lambda_\pm(a)$  given by formula \e{eq:eig}.
 In this case, we  have
 \begin{equation}
 (E_{a}(\{\lambda_{\pm}(a)\})e_0, e_0)= \frac{a^2-2}{2(a^2-1)}.
\label{eq:EiG}\end{equation}
    \end{enumerate}
 \end{theorem}

 Note that the measure $d(E_{a}(\lambda)e_0, e_0)$ is invariant with respect to the reflection $\lambda\mapsto -\lambda$.  Of course for $a=1$, we recover  expression  \e{eq:ZE+} for the spectral measure of the discrete Schr\"odinger operator $H_{+}=H_{1}$ and the corresponding Chebyshev polynomials ${\sf U}_{n} (\lambda)$ of the second kind. If $a=\sqrt{2}$, then \e{eq:ZEa} yields the standard expression  for the orthogonality measure of   Chebyshev polynomials ${\sf T}_{n} (\lambda)$ of the first kind.  According to formula \e{eq:ZEa}  the generalized Chebyshev polynomials $\Ch_{n} (\lambda;a)$ fit into the class of
 Szeg\"o  polynomials (see, for example, \S 10.21   of the book \cite{BE} or \S 2.6 of the book   \cite{Sz})   for all $a\in (0,\sqrt{2}\,]$  (but not for $a>\sqrt{2}$).

 In view of general formula \e{eq:J4}, Theorem~\ref{TZab} implies the next result.
 
  \begin{corollary}\label{TZac}
For all $n,m\in {\Bbb Z}_{+}$, the following representations
\[
d(E_{a}(\lambda)e_n, e_m)= 2 \pi^{-1} \frac{a^2\sqrt{1-\lambda^2 }}{a^4-4 (a^2-1) \lambda^2}\Ch_{n} (\lambda;a)\Ch_m(\lambda;a)d\lambda, \q \lambda\in (-1,1),  
\]
and
 \[
 (E_{a}(\{\lambda_{\pm}(a)\})e_n, e_m)= \frac{a^2-2}{2(a^2-1)}\Ch_{n} (\lambda_{\pm}(a);a)\Ch_m(\lambda_{\pm}(a);a), \q
 a>\sqrt{2},
 \]
 are true.
 \end{corollary}
 
 The following identity is an obvious consequence of Corollary~\ref{TZac}:
 \begin{equation}
 2 \pi^{-1}\int_{-1}^1 \frac{a^2\sqrt{1-\lambda^2 }}{a^4-4 (a^2-1) \lambda^2}\Ch_{n} (\lambda;a)\Ch_m(\lambda;a)d\lambda=
 \d_{n,m}-{\cal I}_{a} ( n, m)
\label{eq:itC}\end{equation}
 where   ${\cal I}_{a} ( n, m)=0$ if $a\leq \sqrt{2}$ and
   \[
{\cal I}_{a} (n, m) = \frac{a^2-2}{2(a^2-1)} \Big(  {\Ch}_{n}(\lambda_+(a);a){\Ch}_m(\lambda_{+} (a);a)  +
 {\Ch}_{n}(\lambda_-(a);a){\Ch}_m(\lambda_- (a);a)  \Big)
\]
if $a >  \sqrt{2}$.

\bigskip

 {\bf 3.4.}
 Following Subsection~2.1 (cf. definition \e{eq:Ua}), we introduce the operators $U_{a}: \ell^2 ({\Bbb Z}_{+})\to L^2 ({\Bbb R}; d\rho_{a})$ by the formula 
$
 (U_{a} e_{n})(\lambda)= \Ch_{n} (\lambda;a),
$
 but it will be convenient to remove the point part of $d\rho_{a}$ (which is non-trivial for $a>\sqrt{2}$ only). Moreover, making a  change of variables we replace $L^2 ((-1,1); d\rho_{a})$ by  $L^2 (-1,1 )$ ($L^2$ with the Lebesgue measure).
 Thus we set
 \[
(V_{a}f)(\lambda)= \sqrt{\frac{2}{\pi} }  \frac{a\sqrt[4]{1-\lambda^2 }}{\sqrt{a^4-4 (a^2-1) \lambda^2}} f(\lambda)
 \]
so that the operator  $V_{a}: L^2 ((-1,1); d\rho_{a})\to L^2 (-1,1)$ is unitary
 and introduce the operator $F_{a}: \ell ^2 ({\Bbb Z}_{+})\to L^2 (-1,1 )$ by the formula
  \begin{equation}
 (F_{a} f)(\lambda)= ( V_{a}\chi_{(-1,1)}U_{a} f) (\lambda), \q \lambda\in (-1,1),
\label{eq:UF}\end{equation}
where $\chi_{(-1,1)}$ is the multiplication operator by the characteristic function of the interval $(-1,1)$.
To put it differently, we set
 \begin{equation}
\psi_{n} (\lambda;a)= \sqrt{\frac{2}{\pi} }  \frac{a\sqrt[4]{1-\lambda^2 }}{\sqrt{a^4-4 (a^2-1) \lambda^2}}   \Ch_n(\lambda;a), \q \lambda\in (-1,1).
\label{eq:ps}\end{equation}
Then
 \begin{equation}
 (F_{a} e_{n})(\lambda)=\psi_{n} (\lambda;a), \q \lambda\in (-1,1).
\label{eq:psX}\end{equation}
 The operator $F_{a}^* :  L^2 (-1,1)\to  \ell^2 ({\Bbb Z}_{+})$ adjoint to $F_{a}$  is given by the formula
 \[
(F_{a}^* g)_{n}=   \int_{-1}^1  \psi_{n} (\lambda;a) g(\lambda) d\lambda,\q n\in {\Bbb Z}_{+}.
\]

 According to \e{eq:UF} it follows from the unitarity of the operator $U_{a}$ that
  \begin{equation}
F_{a} F_{a}^*=I,\q  F_{a}^* F_{a}  =  E_{a}(-1,1) ,
\label{eq:ps3}\end{equation}
where $E_{a}(-1,1)$ is the spectral projection of the operator $H_{a}$ corresponding to the interval $(-1,1)$; of course,
$E_{a}(-1,1)=I$ if $a\leq\sqrt{2}$.

 \bigskip

 {\bf 3.5.}
 Finally, we note that $\omega(z)$ and $D_{a} (z)$ are analytic functions on a two sheets Riemann surface. The second sheet is distinguished by the condition $\sqrt{z^{2}-1}<0$ for $z> 1$. With this convention, formula \e{eq:D} for the perturbation determinant  $D_{a} (z)$ remains true on the second sheet. Its
 zeros   are usually interpreted as resonances (also called anti-bound states). In our case these zeros are simple.
 
 Quite similarly to Proposition~\ref{RZE}, one proves the following result.
 
  \begin{proposition}\label{SH}
  If $a\in (0,1)$, then the operator $H_{a}$ has  two resonances at the points $\pm i 2^{-1} a^2(1-a^2 )^{-1/2}$
lying on the imaginary axis. If $1<a\leq\sqrt{2}$, then the operator $H_{a}$ has   two real  resonances at the points
$\pm  2^{-1} a^2(a^2 -1 )^{-1/2}$.
If $a>  \sqrt{2}$,  then the operator $H_{a}$ does not have resonances $($but it has two eigenvalues$)$.
 \end{proposition}
 
 Compared to the operator $A_{\alpha}$ in the space $L^{2} ({\Bbb R}_{+})$ the picture is of course essentially more complicated. Note that the operator $H_{0}$ has a simple eigenvalue at the point $\lambda=0$ (the corresponding $\omega=\pm i$). As $a$ increases, it splits into two resonances lying on the imaginary axis and tending to $\pm i\infty$ as $a\to 1-0$. If  $a\in (1,\sqrt{2} ]$,  these resonances belong to the real axis and tend from $\pm \infty$   to $\pm 1$   as $a$ increases from $1$ to $\sqrt{2} $. For 
$a>\sqrt{2} $,  the resonances become the eigenvalues.

 \section{Scattering theory}

 {\bf 4.1.} 
 The wave operators $W_{\pm} (H_{a}, H_{1})$   for a pair of self-adjoint operators $H_{1}$, $H_{a}$  are defined as strong limits
   \begin{equation}
  W_{\pm} (H_{a}, H_{1})=\slim_{t\to \pm\infty}e^{iH_{a}t}e^{-iH_{1}t}.
\label{eq:WO}\end{equation}
We refer to the book \cite{Ya} for basic notions of scattering theory. Under the assumption of the existence of limits \e{eq:WO}, $ W_{\pm} (H_{a}, H_{1})$ are isometric operators and enjoy the intertwining property
 $H_{a}W_{\pm} (H_{a}, H_{1})=W_{\pm} (H_{a}, H_{1}) H_{1}$.

 For the operators $H_{1}$, $H_{a}$ defined by formula \e{eq:ZP},
 the perturbation \e{eq:Va}  has finite rank ($V_{a}$ has rank $2$). By the classical Kato theorem, in this case the wave operators $W_{\pm} (H_{a}, H_{1})$ exist.
  Moreover,   the wave operators $W_{\pm} (H_{a}, H_{1})$  are complete, that is, their ranges coincide with the absolutely  continuous subspace of the operator $H_{a}$ (in particular, they are unitary if $a\leq\sqrt{2}$). Therefore the scattering operator
 \begin{equation}
 {\bf S}_{a}=W_{+} (H_{a}, H_{1})^* W_{-} (H_{a}, H_{1})
 \label{eq:WOS}\end{equation}
 is unitary and commutes with $H_{1}$: ${\bf S}_{a}H_{1}=H_{1}{\bf S}_{a}$. All these results will also be obtained in  Subsection~4.3 by a direct method relying on  Theorem~\ref{TZab}.  
  
 To define the corresponding scattering matrix, we use the  diagonalization of   $H_{1}$ by the operator $F_{1}$. 
 Since  the scattering operator ${\bf S}_{a}$ commutes with $H_{1}$ and the operator $H_{1}$ has simple spectrum, we have 
  \begin{equation}
  (F_{1} {\bf S}_{a} f)(\lambda)=S_{a} (\lambda ) (F_{1} f)(\lambda),\q \lambda\in (-1,1),
   \label{eq:SM}\end{equation}  where
$S_{a} (\lambda )\in{\Bbb C}$ and $|S_{a} (\lambda )|=1$. The function $S_{a} (\lambda )$   is known as the scattering matrix for the pair $H_{1}$, $H_{a}$ and the value  $\lambda$ of the spectral parameter. Note that the scattering matrix does not depend on the diagonalization of $H_{1}$ since all its diagonalizations  have the form $\theta(\lambda ) F_{1}$ where $\theta (\lambda )\in{\Bbb C}$ and $|\theta (\lambda )|=1$.

The simplest way to calculate the scattering matrix is to use its abstract expression (this is a particular case of the general Birman-Kre\u{\i}n formula; see, e.g., \S 8.4 of   \cite{Ya}) via the corresponding perturbation determinant $D_{a}(z)$:
 \[
S_{a}(\lambda)= \frac{D_{a}(\lambda- i0)}{D_{a}(\lambda+i0)}, \q \lambda \in (-1,1).
\]
It follows from formula \e{eq:ome1} and  the representation \e{eq:D}  that the limits 
\begin{equation}
D_{a}(\lambda\pm i0)=1 +(1-a^2) (2\lambda^2-1\mp 2i \lambda \sqrt{1-\lambda^2}), \q \lambda \in (-1,1),
\label{eq:D1}\end{equation}
exist, are continuous functions of $\lambda \in (-1,1)$ and $D_{a}(\lambda\pm i0)\neq 0$.
It is also easy to see that $|D_{a}(\lambda\pm i0)|=\sqrt{a^4 + 4 (1-a^2) \lambda^2}$ and
\[
\lim_{\lambda\to\pm 1}D_{a}(\lambda\pm i0)=2-a^2.
\]
 Let us state the result obtained.

 \begin{theorem}\label{S}
The scattering matrix $S_{a}(\lambda)$ for the pair of the  operators $H_{1}$, $H_{a}$ is given by the formula
 \begin{equation}
S_{a}(\lambda)=\frac{1 +(1-a^2) (2\lambda^2-1 +2i \lambda \sqrt{1-\lambda^2})}{1 +(1-a^2) (2\lambda^2-1- 2i \lambda \sqrt{1-\lambda^2})}, \q \lambda \in (-1,1).
\label{eq:SS3}\end{equation} 
In particular,  $S_{a}(-\lambda)=\ov{S_{a}(\lambda)}$ and $S_{a}(0)=1$. Moreover, $S_{a}(\lambda)\to 1$  as $\lambda\to\pm 1$ unless $a=\sqrt{2}$.
 \end{theorem}
 
  \begin{corollary}\label{Sd}
Let $a=\sqrt{2}$. Then, for all $\lambda \in (-1,1)$,
 \[
D_{\sqrt{2}}(\lambda\pm i0)=2\sqrt{1-\lambda^2} (\sqrt{1-\lambda^2}\pm i\lambda)
\]
and
 \[
S_{\sqrt{2}}(\lambda)= \frac{\sqrt{1-\lambda^2}- i\lambda}{\sqrt{1-\lambda^2}+ i\lambda} .
\]
In particular, $S_{a}(\lambda)\to -1$  as $\lambda\to\pm 1$.
 \end{corollary}

\medskip

{\bf 4.2.}
Let us now  discuss the  Kre\u{\i}n spectral shift function $\xi_{a}(\lambda)$ for the pair of the operators $H_{1}$ and $H_{a}$ (see, for example, Chapter~8 of the book \cite{Ya}, for all necessary definitions). In terms of the  perturbation determinant \e{eq:PD}, this function can be introduced  by the formula
 \begin{equation}
\xi_{a}(\lambda)=\pi^{-1}\lim_{\varepsilon\to +0} \arg D_{a} (\lambda+i\varepsilon),\q \lambda\in{\Bbb R}.
\label{eq:SF}\end{equation}
Recall that the  perturbation determinant $D_{a} (z)$ is given by the explicit formula \e{eq:D}.
Since $D_{a}(z)\to 1$ as $|z|\to\infty$ and $D_{a}(z)\neq 0$  for $\Im z\neq 0$, the branch of $\arg D_{a} (z)$  is correctly fixed for $\Im z>0$ by the condition 
$\arg D_{a}(z)\to 0$ as $|z|\to\infty$. In the framework of abstract  scattering theory, the limits in \e{eq:SF} exist for almost all $\lambda\in{\Bbb R}$, but in our particular case they exist for all $\lambda\in{\Bbb R}$ except eventually the points $\pm 1$ and $\lambda_\pm(a)$.
If $a\leq\sqrt{2}$, then $\xi(\lambda)=0$ for all $|\lambda| >1$. In the case $a>\sqrt{2}$, the function  $\xi_{a}(\lambda)=0$   for $|\lambda|> |\lambda_{\pm}(a)|$, $\xi_{a}(\lambda)=1$ for $\lambda\in (1, \lambda_+(a))$ and $\xi_{a}(\lambda)=-1$ for $\lambda\in (\lambda_{-}(a), -1)$.
 
 Let us now consider the function $\xi_{a}(\lambda)$ for $\lambda\in (-1,1)$. Since   $D_{a}(\lambda + i0)$ is a continuous   function of $\lambda \in (-1,1)$ and $D_{a}(\lambda + i0)\neq 0$, the spectral shift function $\xi_{a}(\lambda)$ depends also continuously on  $\lambda \in (-1,1)$. By \e{eq:id2},  \e{eq:D}, $\Im  D_{a}(iy)= 0$
  and thus $D_{a}(iy)>0$ for all $y\geq 0$ so that  $\xi_{a}(0)  =0$.  It also follows   from \e{eq:D1}  that $D_{a}(-\lambda+i0)= \ov{D_{a} (\lambda+i0)}$ whence $\xi_{a}(-\lambda)=-\xi_{a}(\lambda)$. So it suffices to study $\xi_{a}(\lambda)$ for $\lambda\in ( 0,1)$.

 Putting together formulas \e{eq:D1} and \e{eq:SF}, we see that
  \begin{equation}
  \tan (\pi \xi_{a}(\lambda))= \frac{2(a^2-1)\lambda \sqrt{1-\lambda^2}}{1+ (1-a^2) (2\lambda^2-1)} .
\label{eq:SF1}\end{equation}
Let us now consider separately the cases $a<\sqrt{2}$,  $a>\sqrt{2}$ and  $a=\sqrt{2}$. In the first case the denominator in
\e{eq:SF1} does not equal zero for $\lambda\in(0,1)$, and  hence formula \e{eq:SF1} can be rewritten as
 \begin{equation}
  \xi_{a}(\lambda)=  \frac{1}{\pi}\arctan \frac{2(a^2-1)\lambda \sqrt{1-\lambda^2}}{1+ (1-a^2) (2\lambda^2-1)} .
\label{eq:SF2}\end{equation}
Obviously, $  \xi_{a}(\lambda)\to 0$ as $\lambda\to  1$. It is also easy to see  that $\xi_{a} (\lambda)<0$ for $a\in(0,1)$
and $\xi_{a} (\lambda)> 0$ for $a\in( 1, \sqrt{2})$. Moreover,  $\xi_{a} (\lambda)$ 
has one extremum at the point $\lambda=a/\sqrt{2}$ where
\[
\xi_{a} (a/\sqrt{2})= \frac{1}{\pi}\arctan \frac{(a^2-1)   }{ a\sqrt{2-a^2}}.
\]
This point is the
minimum of $\xi_{a} (\lambda)$ for $a\in(0,1)$ and its maximum for $a\in( 1, \sqrt{2})$.

If $a>\sqrt{2}$, then the denominator in
\e{eq:SF1}   equals zero at the point $a 2^{-1/2} (a^{2}-1)^{-1/2}\in (0,1)$. In this case
 it follows from \e{eq:SF1}   that
 \begin{equation}
  \xi_{a}(\lambda)= \frac{1}{\pi}\arccot \frac{1- (a^2-1) (2\lambda^2-1)} {2(a^2-1)\lambda \sqrt{1-\lambda^2}}  .
\label{eq:SF3}\end{equation}
In particular,  we see that $ \xi_{a}'(\lambda)> 0$ and $ \xi_{a}(\lambda)\to 1$ as $\lambda\to 1$.

Let us summarize the results obtained.

 \begin{theorem}\label{SSF}
Let the Jacobi operator $H_{a}$ in the space $\ell^2 ({\Bbb Z}_{+})$ be defined by formula \e{eq:ZP+}, and let $\xi_{a} (\lambda)$ be the  spectral shift  function for the pair $H_{1}, H_{a}$.
 Then $\xi_{a}(\lambda)$ is an odd function of $\lambda\in{\Bbb R}$ and the following results are true.
 
  \begin{enumerate}[\rm(i)]

\item
Let $a\in (0,\sqrt{2})$. Then $\xi_{a} (\lambda)$ is given by formula \e{eq:SF2} for $\lambda\in [0,1)$ and 
$\xi_{a} (\lambda)=0$ for $\lambda \geq 1$.

\item
Let $a> \sqrt{2}$. Then $\xi_{a} (\lambda)$ is given by formula \e{eq:SF3} for $\lambda\in [0 ,1)$,
$\xi_{a} (\lambda)=  1$ for $\lambda\in [  1, \lambda_{+} (a))$ and $\xi_{a} (\lambda)=0$ for $\lambda \geq \lambda_{+} (a)$.

\item
In the intermediary case $a =\sqrt{2}$, we have
\[
  \xi_{\sqrt{2}}(\lambda)=  \frac{1}{\pi}\arctan \frac{ \lambda }{ \sqrt{1-\lambda^2}} , \q \lambda\in [0 ,1),
\]
$\xi_{\sqrt{2}}(\lambda)\to1/2$ as $\lambda\to 1$
and $  \xi_{\sqrt{2}}(\lambda)=  0$ for $\lambda>1$. 

    \end{enumerate}
 \end{theorem}
 
Note that
the  perturbation determinant \e{eq:PD} can be recovered from the spectral shift function     by the formula
 \[
 \ln D_{a} (z)=\int_{-\infty}^\infty \frac{\xi_{a}(\lambda)}{\lambda-z}d\lambda .
\]
Substituting here expression \e{eq:D} for $D_{a} (z)$and the expressions of Theorem~\ref{SSF} for $\xi_{a}(\lambda)$, we obtain identities parametrized by $a\in {\Bbb R}_{+}$ which do not look quite obvious.

\medskip

 {\bf 4.3.}  
The wave operators and the scattering matrix can be expressed via eigenfunctions (of the continuous spectrum) of the operator $H_{a}$. These eigenfunctions can be constructed
in terms of the generalized Chebyshev polynomials.
It is convenient to introduce  the operator  ${\sf A}$ of multiplication by $\lambda$ in the space $L^2 (-1,1)$.
Let the operators $F_{a}$ be defined by formula \e{eq:UF}.
According to the  equation \e{eq:IP} the intertwining property
  \begin{equation}
H_{a}F_{a}^*=F_{a}^* {\sf A}
\label{eq:inw1}\end{equation}
holds.

 Next, we find a relation between the operators $F_{a}$ and the wave operators $W_\pm (H_{a}, H_{1})$. We use  the following elementary result.
 
  \begin{lemma}\label{CheW}
  Let  the function $\omega (\lambda\pm i0)$ be given by formula \e{eq:ome1}, and let $g\in C_{0}^\infty (-1,1)$. Then, for an arbitrary $p\in {\Bbb Z}_{+}$ and $t\to \mp \infty$, we have 
   \begin{equation}
\big| \int_{-1}^1 \omega (\lambda\pm i0)^n e^{-i\lambda t} g(\lambda) d\lambda \big|
\leq C_{p} (n+ |t|)^{-p},\q n\in {\Bbb Z}_{+},
\label{eq:opm1}\end{equation}
with some constant  $C_{p}$ depending on $p$ only.
  \end{lemma}

  \begin{pf}
  Let us integrate by parts in \e{eq:opm1} and
  observe that
  \[
    \frac{d}{d\lambda} \big(n\ln  \omega (\lambda\pm i0)  -i\lambda t)
  =\pm i \big( \frac{n} {\sqrt{1-\lambda^2}} \mp t \big) .
  \]
  If $\mp t >0$, the modulus of this expression is bounded from below by $n+|t|$ which yields estimate \e{eq:opm1} for $p=1$. Similarly,  integrating in \e{eq:opm1} $p$ times by parts, we obtain   estimate  \e{eq:opm1} for the same value of $p$.
      \end{pf}
      
    Let us now set  
     \begin{equation}
 \sigma_\pm(\lambda;a) = \frac{a^2+2 (1- a^2) \lambda^2\mp
 i2 (a^2-1)\lambda \sqrt{1-\lambda^2}} {\sqrt{a^4-4 (a^2-1)\lambda^2}}  .  
 \label{eq:CH1r}\end{equation} 
 Clearly, $ |\sigma_\pm(\lambda;a)|=1$ so that the operator $ \Sigma_\pm(\lambda;a)$ of multiplication by $ \sigma_\pm(\lambda;a)$ is unitary in the space $ L^2 (-1,1)$.

      \begin{lemma}\label{CheW1}
 For all $g\in L^2 (-1,1)$, we have
    \begin{equation}
 \lim_{t \to\pm\infty} \| (F_{a}^* \Sigma_{\pm} (a) -F_{1}^*)e^{-i{\sf A}t}g\|=0
 \label{eq:CH1rr}\end{equation}
  \end{lemma}
  
    \begin{pf}
It suffices to check  \e{eq:CH1rr} for  $g\in C_{0}^\infty (-1,1)$.
    Let us proceed from definition \e{eq:ps} where  the polynomial ${\Ch}_{n}(\lambda;a)$  is given by formula \e{eq:CH2}
    with  $z=\lambda+i0$ and the numbers $\omega(\lambda\pm i0)$ and $\gamma_{\pm}(\lambda+i0 ;a)$ are given by formulas \e{eq:ome1} and \e{eq:CH12}, respectively:
    \begin{equation}
\psi_{n} (\lambda;a)= \sqrt{\frac{2}{\pi} }  \frac{a\sqrt[4]{1-\lambda^2 }}{\sqrt{a^4-4 (a^2-1) \lambda^2}}   \big( \gamma_{+}(\lambda+i0 ;a)  \omega(\lambda+i0)^{n} +   \gamma_{-}(\lambda+i0 ;a) \omega(\lambda-i0)^{n}\big).
\label{eq:PS}\end{equation} 
In particular, for $a=1$, we have
 \begin{equation}
\psi_{n} (\lambda;1)= \sqrt{\frac{2}{\pi} }  \sqrt[4]{1-\lambda^2 }   \big( \gamma_{+}(\lambda+i0 ;1)  \omega(\lambda+i0)^{n} +  \gamma_{-}(\lambda+ i0 ;1) \,  \omega(\lambda-i0)^{n}\big).
\label{eq:PS1}\end{equation} 
Lemma~\ref{CheW} implies that relation \e{eq:CH1rr} is satisfied if  the coefficient at $ \omega(\lambda\pm i0)^{n}$ in the expression
 \[
\psi_{n} (\lambda;a) \sigma_\pm(\lambda;a)-\psi_{n} (\lambda;1)
\]
is zero. According to \e{eq:PS}, \e{eq:PS1}  this yields the equation
 \[
\frac{a }{\sqrt{a^4-4 (a^2-1) \lambda^2}}    \gamma_{\pm}(\lambda + i0 ;a)  \sigma_\pm(\lambda;a)=  \gamma_{\pm}(\lambda+ i0 ;1) 
\]
for $\sigma_\pm(\lambda;a)$.  Its solution is given by formula \e{eq:CH1r}.
  \end{pf}

   \begin{proposition}\label{CheW2}
 For all $a>0$, the strong limits \e{eq:WO}  exist and
     \begin{equation}
W_{\pm}(H_{a},H_{1})= F_{a}^* \Sigma_{\pm} (a)  F_{1}.
 \label{eq:CH1pr}\end{equation}
  \end{proposition}
  
     \begin{pf}
   We have to check that
          \[
 \lim_{t \to\pm\infty} \| e^{iH_{a}t} e^{-iH_1 t}f - F_{a}^* \Sigma_{\pm} (a) F_{1}f\|=0
\]
 for all $f\in\ell^2({\Bbb Z}_{+})$. In view of the intertwining property \e{eq:inw1} this relation can be rewritten as
   \begin{equation}
 \lim_{t \to\pm\infty} \| e^{-iH_1 t}f - F_{a}^* \Sigma_{\pm} (a)  e^{-i{\sf A}t} F_{1}f\|=0.
 \label{eq:wo1}\end{equation}
 Set now $g=F_{1}f$. Then again in view of the intertwining property \e{eq:inw1} for   $a=1$, we see that relations
 \e{eq:CH1rr} and \e{eq:wo1}  are equivalent.
          \end{pf}

  It follows from relations \e{eq:ps3} and  \e{eq:CH1pr}  that the scattering operator \e{eq:WOS} is given by the equality
  \[
  {\bf S}_{a}= F_{1}^* \Sigma_{+} (a)^{*} \Sigma_{-} (a) F_{1}.
  \]
  Putting together this relation with the definition  \e{eq:SM} of  the scattering matrix, we see that
   \[
S_{a} (\lambda)=  \frac{\sigma_{-}(\lambda ;a)}{\sigma_{+}(\lambda ;a)} .
\]
Thus according to formula \e{eq:CH1r} for $\sigma_\pm (\lambda ;a)$, we recover the representation   \e{eq:SS3}
for $S_{a} (\lambda)$. Still another proof of   \e{eq:SS3} will be given in Subsection~6.3.

\section{Trace identities and moment problems}

{\bf 5.1.} 
Let us first find matrix elements $ ( R_{a} (z)  e_n,e_m)$ for all $n,m\in{\Bbb Z}_{+}$. Recall that for $a=1$ they are given by formula \e{eq:Z4+} and denoted $r_{n,m} (z)$. As usual, $\omega(z)$ and $D_{a}(z)$ are the functions \e{eq:ome} and \e{eq:D}, respectively.

 \begin{proposition}\label{TTnm}
   Let     $n\geq 1$. Then
   \begin{equation}
( R_{a} (z)  e_n,e_m)=r_{n,m} (z) -4 (a^2-1)  D_{a}(z)^{-1} z \omega(z)^{n+m+2} 
\label{eq:TFnm}\end{equation}
for $m\geq 1$ and   
   \begin{equation}
( R_{a} (z)  e_n,e_{0})= r_{n,0} (z) -2 (a-1)  D_{a}(z)^{-1}\omega(z)^{n+2} (2z+ a \omega(z)).
\label{eq:TFn0}\end{equation}
 \end{proposition}
 
  \begin{pf}
We  follow the   scheme of the proof of Theorem~\ref{TTab}.   If $ n\geq 1$ and $m\geq 1$, then  Lemma~\ref{RZr} and  Proposition~\ref{TZa}  imply that
\begin{multline}
D_{a}(z)(R_{1}(z)  T_{a} (z)R_{1}(z) e_n, e_{m})
\\
 = 4 D_{a}(z)\omega(z)^{n+m+2}  (  T_{a} (z) (e_0+2z  e_{1}), e_0+2 \bar{z}   e_{1})
 \\
=  2 \omega(z)^{n+m+2}  ( t_{0,0}(z)+4z t_{0,1}(z)+4z^2  t_{1,1}(z)).
\label{eq:T7Aa}\end{multline}
Substituting here expressions \e{eq:Ta7}, we obtain that
   \begin{equation}
D_{a}(z)(R_{1}(z)  T_{a} (z)R_{1}(z) e_n, e_m) = 4 (a^2-1)  z  \omega(z)^{n+m+2}.
\label{eq:T7Aaa}\end{equation}

  If $n\geq 1$, but $m=0$, then instead of \e{eq:T7Aa} we have the formula
  \begin{multline*}
D_{a}(z)(R_{1}(z)  T_{a} (z)R_{1}(z) e_n, e_0)
\\
 = 4 D_{a}(z)\omega(z)^{n+ 2}  (  T_{a} (z) (e_0+2z  e_{1}), e_0+\omega ( \bar{z} )  e_{1})
 \\
=  2 \omega(z)^{n+ 2}  ( t_{0,0}(z)+2z t_{1,0}(z)+\omega(z ) t_{0,1}(z)+ 2z \omega(z)  t_{1,1}(z)).
\end{multline*}
  Substituting here expressions \e{eq:Ta7}, we obtain that
   \begin{equation}
D_{a}(z)(R_{1}(z)  T_{a} (z)R_{1}(z) e_n, e_0) = 2 (a-1)    \omega(z)^{n+2} (2z+  a \omega(z)).
\label{eq:T7Abb}\end{equation}

In view of formula \e{eq:T7A} relations  \e{eq:TFnm} and \e{eq:TFn0} follow from \e{eq:T7Aaa} and \e{eq:T7Abb}, respectively.
  \end{pf}
 
   \begin{corollary}\label{TTnx}
   Let     $n\geq 1$. Then
   \[
( R_{\sqrt{2}} (z)  e_n,e_m)=r_{n,m} (z) - 2 (z^{2}-1)^{-1/2} z \omega(z)^{n+m+1} 
\]
for $m\geq 1$ and
\[
( R_{\sqrt{2}} (z)  e_n,e_{0})= r_{n,0} (z) -(2-\sqrt{2}) (z^{2}-1)^{-1/2}   \omega(z)^{n+1} (\sqrt{2} z+  \omega(z)).
\]
  \end{corollary}
 
 \medskip
 
{\bf 5.2.} 
The following result is a direct consequence of relation \e{eq:J6}  and  Theorem~\ref{TZab}.

\begin{proposition}\label{Che}
   Let the coefficients $( R_{a} (z)  e_n,e_m)$ be defined by formulas \e{eq:Ta7X}, \e{eq:TFnm} and \e{eq:TFn0}, and let the numbers $\lambda_{\pm} (a)$ be given by equality  \e{eq:eig}.  Then
   for all $a> 0$ and all $n,m\in{\Bbb Z}_{+}$, the generalized Chebyshev polynomials
${\Ch}_{n}(\lambda;a)$ satisfy the identity
    \begin{equation}
\frac{2   a^2 } {\pi}  \int_{-1}^1\frac{ {\Ch}_{n}(\lambda;a){\Ch}_m(\lambda;a)\sqrt{1-\lambda^2 }}{(\lambda-z) (a^4-4 (a^2-1) \lambda^2)}d\lambda  =( R_{a} (z)  e_n,e_m) -I_{a} (z; n, m)
\label{eq:itC1}\end{equation}
where $I_{a} (z; n, m)=0$ if $a\leq \sqrt{2}$ and
   \[
I_{a} (z; n, m) = \frac{a^2-2}{2(a^2-1)} \Big(\frac{ {\Ch}_{n}(\lambda_+(a);a){\Ch}_m(\lambda_{+} (a);a) }{\lambda_{+} (a) -z}+
\frac{ {\Ch}_{n}(\lambda_-(a);a){\Ch}_m(\lambda_- (a);a) }{\lambda_-(a) -z}\Big)
\]
if $a >  \sqrt{2}$.
 \end{proposition}

 Let us  state a consequence of this result  for $a=1$ and $a=\sqrt{2}$ when $ {\Ch}_{n}(\lambda;1) ={\sf U}_{n}(\lambda)$ and  $ {\Ch}_{n}(\lambda;\sqrt{2}) ={\sf T}_{n}(\lambda)$
 are the classical  Chebyshev polynomials. The first assertion follows, actually, from Proposition~\ref{RZ+}.
  
  \begin{corollary}\label{CheU}
  For all $n,m\geq 0$, we have
  \[
 \int_{-1}^1 (\lambda-z)^{-1}{\sf U}_{n}  (\lambda){\sf U}_m (\lambda) \sqrt{1-\lambda^2 }d\lambda
= \pi \frac{(z-\sqrt{z^2 -1})^{n+m+2}-(z-\sqrt{z^2 -1})^{|n-m|} }{2\sqrt{z^2 -1}} .
\]
 \end{corollary}
 
 The second assertion  requires Corollary~\ref{TTnx} only.

 \begin{corollary}\label{CheT}
 The integral
  \[
\frac{1} {\pi}  \int_{-1}^1\frac{ {\sf T}_{n}(\lambda){\sf T}_{m}(\lambda) }{(\lambda-z) \sqrt{1-\lambda^2 }}d\lambda  
\]
equals $-(z^2-1)^{-1/2}$ if $n=m=0$. It equals $-  \sqrt{2}\omega(z)^{n }  $
if $n\geq 1$ but $m=0$, and it equals
\[
r_{n,m}(z) -\frac{2  \omega(z)^{n+m+1}}{\sqrt{z^2-1}}
\]
if $n\geq 1$, $m\geq 1$.
 \end{corollary}
 
 The identities of Corollaries~\ref{CheU} and \ref{CheT} are probably well-known although we were unable to find them in the literature. 
 
 Of course taking the limit $z\to\infty$ in \e{eq:itC1} and equating   the coefficients at $z^{-1}$, we recover identity \e{eq:itC}.
 
  \medskip
 
{\bf 5.3.} 
To calculate  the trace
  \[
\tr \big(R_{a} (z) - R_{1} (z)\big)=\sum_{n=0}^\infty \big((R_{a} (z) e_{n}, e_{n})- (R_1(z) e_{n}, e_{n})\big),
\]
it is convenient to use a link between the trace and the perturbation determinant \e{eq:PD}:
 \begin{equation}
\tr \big(R_{a} (z) - R_{1} (z)\big)=-\frac{D_{a}'(z)}{D_{a}(z)}.
\label{eq:trD}\end{equation}
Differentiating  expression \e{eq:D}, we see that
 \[
 D_{a}'(z)=2 (1-a^2) \omega(z) \omega'(z)
\]
where, by definition \e{eq:ome}, 
\[
\omega'(z)=-\frac{\omega(z)}{\sqrt{z^2-1}}.
\]
Substituting these expressions into \e{eq:trD}, we obtain the following result.

\begin{proposition}\label{tr}
For all $z\in{\Bbb C}\setminus [-1,1]$, we have
   \begin{equation}
\tr \big(R_{a} (z) - R_{1} (z)\big)=\frac{2(1-a^2) \omega (z)^2}{D_{a}(z) \sqrt{z^2-1}}.
\label{eq:tr1}\end{equation}
 \end{proposition}
 
 It is now easy to obtain a generating function for the sequence $\tr \big(H_{a}^n - H_{1}^n\big)$. Indeed, the left-hand side
 of \e{eq:tr1} has 
 the asymptotic expansion as $z\to\infty$ (for example, along the positive half-line):
   \[
\tr \big(R_{a} (z) - R_{1} (z)\big)=-\sum_{n=1}^\infty z^{-n-1}\tr \big(H_{a}^n - H_{1}^n\big).
\]
Putting now $\zeta=z^{-1}$ and accepting that the function $\sqrt{1-\zeta^2}=1$ at $\zeta =0$, we can state a consequence of Proposition~\ref{tr}.

\begin{corollary}\label{trc}
We have
   \[
 \sum_{n=1}^\infty \zeta^{n}\tr \big(H_{a}^n - H_{1}^n\big)=\frac{2(a^2-1)}{\sqrt{1-\zeta^2}}
 \frac{\big( 1-\sqrt{1-\zeta^2}\big)^2 }{\zeta^2   + (1-a^2)\big( 1-\sqrt{1-\zeta^2}\big)^2}
\]
where the series in the left-hand side converges for $|\zeta|<1$ if $a\leq\sqrt{2}$ and  for $|\zeta|<2 a^{-2} \sqrt{a^2-1}$ if $a>\sqrt{2}$. In particular,   $\tr \big(H_{a}^n- H_{1}^n\big)=0$ if $n$ is odd. 
 \end{corollary}
 
 We note the paper \cite{Case} where the expressions for $\tr \big(H_{a}^n - H_{1}^n\big)$ were obtained for general Jacobi operators $H$. Of course these expressions are less explicit than for the  particular case of the operators $H=H_{a}$.
 
   \medskip
 
{\bf 5.4.} 
Finally, we consider the moments
   \begin{equation}
\kappa_{n}(a)=\int_{-\infty}^\infty \lambda^n d\rho_{a} (\lambda)
\label{eq:MM}\end{equation}
of the measure $d\rho_{a} (\lambda)= d(E_a (\lambda) e_{0}, e_{0})$. By Theorem~\ref{TZab}, this measure is absolutely continuous on the interval $[-1,1]$ and has two atoms at the points $\lambda_{\pm} (a)$ if $a>\sqrt{2}$. We recall that the points
$\lambda_{\pm} (a)$    and the measures $\rho_{a} (\{\lambda_{\pm} (a)\}) $ are given by formulas \e{eq:eig} and \e{eq:EiG}, respectively. Of course $\kappa_{n}(a)=0$ for $n$ odd because the measure  $d\rho_{a} (\lambda)$ is even.

We will find an explicit expression for a generating function of the sequence $\kappa_{n}(a)$.   It follows from Theorem~\ref{TTab} that
\[
\int_{-\infty}^\infty (\lambda-z)^{-1}d\rho_{a} (\lambda)= \frac{2 }{( a^2-2)z- a^2\sqrt{z^2-1}}.
\]
Using the same arguments as those  used for the proof of
Corollary~\ref{trc}, we arrive at the next result.

\begin{proposition}\label{PM}
We have
   \[
 \sum_{n=0}^\infty \kappa_{n}(a) \zeta^{n} =\frac{2 } {a^2\sqrt{1-\zeta^2} +2-a^2 }
\]
where the series in the left-hand side converges for $|\zeta|<1$ if $a\leq\sqrt{2}$ and  for $|\zeta|<2 a^{-2} \sqrt{a^2-1}$ if $a>\sqrt{2}$.  
 \end{proposition}

 Of course, Proposition~\ref{PM} is consistent with   the well-known expressions
 \[
\kappa_{2m}(1)= \frac{(2m-1)!!} {(m+1)! \, 2^m}, \q  \kappa_{2m}(\sqrt{2})= \frac{(2m-1)!!} {m ! \, 2^m}
  \]
  for the classical Chebyshev polynomials.
  It follows from the Stirling formula (see, e.g., formula (1.18.4) in \cite{BE})  that
     \begin{equation}
\kappa_{2m}(1)= \frac{1} {\sqrt{\pi}m^{3/2}}\big( 1+O ( \frac{1} {m })\big),\q  \kappa_{2m}(\sqrt{2})= \frac{1} {\sqrt{\pi}m^{1/2}}\big( 1+O ( \frac{1} {m })\big)
\label{eq:MM3}\end{equation}
as $m\to\infty$. Since the asymptotics   as $m\to \infty$ of the integral \e{eq:MM} is determined by  neighborhoods of the points $\lambda=\pm 1$, it is easy to deduce from expression \e{eq:ZEa} for $d\rho_{a}(\lambda)$ and the first relation \e{eq:MM3} that for all $a <\sqrt{2}$
 \begin{equation}
\kappa_{2m}(a)= \big( \frac{a}{a^2-2} \big)^2\frac{1} {\sqrt{\pi}m^{3/2}}\big( 1+O ( \frac{1} {m })\big),\q   m\to\infty.
\label{eq:MM4}\end{equation}
If $a >\sqrt{2}$, then relations \e{eq:eig} and \e{eq:EiG}   imply that $\kappa_{2m}(a)$ tend to infinity exponentially:
\[
\kappa_{2m}(a)= \frac{a^2-2}{a^2-1}  \Big( \frac{a^4}{4 (a^2-1)}  \Big)^m
  + O(1) ,\q   m\to\infty.
\]
  
   \medskip
   
   \section{Miscellaneous}
 
{\bf 6.1.} 
Let us now briefly discuss the connection of the moment problem with Hankel operators. We denote by ${\cal D}\subset
\ell^2 ({\Bbb Z}_{+})$ the set of vectors $f=(f_{0}, f_{1},\ldots)$ with only a finite number of non-zero components $f_{n}$. The Hankel operator $K$ corresponding to a sequence $\kappa_{n}$ is formally defined   by the relation
 \[
  (K f)_{n}=\sum_{m =0}^\infty \kappa_{n+m}f_{m}, \q \forall f\in{\cal D}.
\]
The classical Hamburger theorem (see, e.g.,  the book \cite{AKH}) states that the Hankel quadratic form 
 \[
\kappa[f,f]:=  \sum_{n,m=0}^\infty \kappa_{n+m}f_{m}\ov{f_{n}}\geq 0, \q \forall f\in{\cal D},
\]
if and only if
   \begin{equation}
\kappa_{n} =\int_{-\infty}^\infty \lambda^n d \rho (\lambda)
\label{eq:Ha2}\end{equation}
with some nonnegative measure $d \rho (\lambda)$. 

If the form $\kappa[f,f]$ is closable in   $\ell^2 ({\Bbb Z}_{+})$, then using the Friedrichs construction (see, e.g., the book \cite{BS})  one can standardly associate 
with it  a nonnegative operator $K=K(\kappa)$  in this space. The necessary and sufficient condition obtained in  \cite{Yunb} for the form determined by  sequence 
\e{eq:Ha2} to be closable is that $\supp \rho\subset [-1,1]$ and $\rho(\{1\})= \rho(\{-1\})=0$. For the sequence $\kappa_{n} (a)$ defined by \e{eq:Ha2} with the measure $d \rho_{a} (\lambda)$ given by  \e{eq:ZEa}, this condition is satisfied if and only if $a\leq\sqrt{2}$.

Let $\supp \rho\subset [-1,1]$  and $\rho(\{1\})= \rho(\{-1\})=0$. Then the Widom theorem \cite{Widom} states that the operator $K$ with matrix elements \e{eq:Ha2} is bounded (resp. compact) if and only if $\rho(1-\varepsilon,1)= O(\varepsilon)$ and 
$\rho(-1, -1-\varepsilon)= O(\varepsilon)$ (resp, $\rho(1-\varepsilon,1)= o(\varepsilon)$ and 
$\rho(-1, -1-\varepsilon)= o(\varepsilon)$) as $\varepsilon\to 0$. Therefore the Hankel operator $K (a) $ with the matrix elements $\kappa_{n} (a)$ is unbounded if $a=\sqrt{2}$.  On the contrary, the operators $K (a) $ are compact if $a<\sqrt{2}$.

\medskip

 {\bf 6.2.} 
 Let us now discuss the behavior of the spectral measure $d\rho_{a} (\lambda)$, of the scattering matrix $S_{a} (\lambda)$  and of the spectral shift function $\xi_{a} (\lambda)$ in the   limits $a\to \infty$ and  $a\to 0$. Both these limits are very singular.

 If  $a\to \infty$ (the   ``large coupling" limit), then according to
 \e{eq:EiG}
  \[
\rho_{a}(\{\lambda_{\pm}(a)\})= \frac{a^2-2}{2(a^2-1)}\to 1/2
\]
and therefore $\rho_{a}(-1,1) \to 0$. On the contrary, according to \e{eq:SS3} the  scattering matrix $S_{a} (\lambda)$ defined for $\lambda\in (-1,1)$ has a non-trivial limit:
 \[
\lim_{a\to\infty}S_{a}(\lambda) = \frac{ 2\lambda^2-1 +2i \lambda \sqrt{1-\lambda^2}} { 2\lambda^2-1 -2i \lambda \sqrt{1-\lambda^2}}\neq 1
\]
(if $\lambda\neq 1$).
Similarly, according to \e{eq:SF3}
 \[
\lim_{a\to\infty}\xi_{a}(\lambda) =  \frac{  1}{\pi}\arccot \frac{1-  2\lambda^2 } {2 \lambda \sqrt{1-\lambda^2}} ,\q \lambda\in (0,1).
\]

Let now $a\to 0$. By \e{eq:ZP+},  the operator $H_0$ can be considered as   the orthogonal sum   of the operators acting in the space ${\Bbb C}\oplus \ell^2({\Bbb N})$ where ${\Bbb N}=(1,2,\ldots)$. The restriction of $H_{0}$ on ${\Bbb C}$ is zero, and its restriction on $\ell^2({\Bbb N})$ is again given by matrix \e{eq:ZP+} where $a=1$.
According to 
 \e{eq:ZEa} for every $\varepsilon\in (0,1)$  we have
 \[
 \rho_{a}(-1,-1+\varepsilon) \to 0, \q  \rho_{a}( 1-\varepsilon,1) \to 0
 \]
  and therefore $\rho_{a}(-\varepsilon,\varepsilon) \to 1$ as $a\to 0$.  It follows from \e{eq:SS3} that the  scattering matrix $S_{a} (\lambda)$  has a non-trivial limit:
 \[
\lim_{a\to 0} S_{a}(\lambda) = \frac{ \lambda +i   \sqrt{1-\lambda^2}} { \lambda - i   \sqrt{1-\lambda^2}}\neq 1 ,\q \lambda\in (-1,1).
\]
Similarly, it follows from \e{eq:SF2} that
 \[
\lim_{a\to 0}\xi_{a}(\lambda) = 
 -\frac{1}{\pi}\arctan \frac{ \sqrt{1-\lambda^2}}{\lambda}  ,\q \lambda\in (-1,1).
\]

\medskip

 {\bf 6.3.} 
 The  expression for
the scattering matrix can also be obtained in terms of the operator $T_{a}(z)$ defined by formula \e{eq:Ta}.  Let ${\bf t}_{a}(\lambda,\mu;z)$ be the integral kernel of the operator  $F_{1} T_{a}(z)F_{1}^*$ where $F_{1}$ is given by \e{eq:psX}
with 
\[
\psi_{n}(\lambda;1)=\sqrt{ \frac{ 2}{\pi}}\sqrt[4]{1-\lambda^{4}}{\sf U}_n(\lambda)
\]
(recall that ${\sf U}_n(\lambda)={\Ch}_{n} (\lambda;1)$).
 Then 
 \begin{equation}
S_{a}(\lambda)=1-2\pi i \, {\bf t}_{a}(\lambda,\lambda;\lambda+ i0),\q \lambda\in (-1,1).
\label{eq:S}\end{equation} 
This formula has been obtained in the paper \cite{Fa} in the framework of the Friedrichs-Faddeev model and discussed in the book \cite{Ya} in   a more general setting.

It follows from formula \e{eq:Ta6} that
\begin{multline*}
\pi D_{a}(z) {\bf t}_{a}(\lambda,\mu;z) = (1-\lambda^2)^{1/4} (1-\mu^2)^{1/4}  
\Big(t_{0,0} (z) {\sf U}_{0}(\lambda)  {\sf U}_{0}(\mu)+ t_{0,1} (z) {\sf U}_{0}(\lambda)  {\sf U}_{1}(\mu)
\\
+   (t_{1,0} (z) {\sf U}_{1}(\lambda)  {\sf U}_{0}(\mu)+ t_{1,1} (z) {\sf U}_{1}(\lambda)  {\sf U}_{1}(\mu)\Big).
\end{multline*}
Since $ {\sf U}_{0}(\lambda) =1$, $ {\sf U}_1(\lambda) =2\lambda$, representation   \e{eq:S} implies that
 \begin{equation}
S_{a}(\lambda)=1-\frac{2i\sqrt{1-\lambda^2}}{D_{a}(\lambda+i0)}\big(t_{0,0} (\lambda+i0)  + 2(t_{0,1} (\lambda+i0) +   t_{1,0} (\lambda+i0))\lambda + 4  t_{1,1} (\lambda+ i0) \lambda^2\big). 
\label{eq:S2}\end{equation} 
Let us now use formulas \e{eq:Ta7} and the second  identity \e{eq:id2}: 
 \begin{multline*}
 t_{0,0} (\lambda+i0)  + 2(t_{0,1} (\lambda+i0) +   t_{1,0} (\lambda+i0))\lambda + 4  t_{1,1} (\lambda+ i0) \lambda^2
\\
 = 2 (a-1) \lambda \big(2-  (a-1) \omega(\lambda+i0)^2 + 2  (a-1) \omega(\lambda+i0)\lambda \big)=2 (a^2-1) \lambda . 
\end{multline*}
Substituting this expression into formula \e{eq:S2}, we recover expression \e{eq:SS3}.

    \medskip
    
     \section{More general perturbations}
   
   
    {\bf 7.1.} 
    Here we briefly consider arbitrary finite rank perturbations $V$ of the free Jacobi operator $H_{1}$. 
    In principle, the scheme of the main part of the paper can be   adjusted to handle  the operator    $H=H_{1}+V$, 
    but here we use another, also quite standard, approach relying on a study of the so called Jost function. 
  
    First we recall its  definition. Let us consider a general Jacobi operator \e{eq:J};  we assume that $a_{n}=1/2$   and  $b_{n}= 0$ for $n\geq N$. In this case, for all $ z\in {\Bbb C}\setminus [-1,1]$ and $n\geq N+1$,  the equations
\begin{equation}
  a_{n-1} u_{n-1}(z) = (z-b_{n} ) u_{n}(z) - a_{n} u_{n+1}(z)  
\label{eq:ab}\end{equation}
are satisfied if $u_m (z)=\omega (z)^m$ for $m\geq N$. Here $\omega(z)$ is defined as usual by \e{eq:ome} and hence $2z= \omega(z)+\omega(z)^{-1}$. Then equations \e{eq:ab} for $n=N $, $n=N-1,\ldots, n=1$,   define recurrently all $u_{N-1} (z)$, $u_{N-1} (z),\ldots, u_{0} (z)$. Since $|\omega(z)|< 1$, the vector $u(z)= (u_{0}(z),u_{1}(z),\ldots)\in \ell^2 ({\Bbb Z}_{+})$. It is called the Jost solution of the equation $H u(z)
=z u(z)$. Finally, for $n=0$,  equation \e{eq:ab}  determines $a_{-1} u_{-1} (z)$; for definiteness,   we put  $a_{-1}=1$. The function $u_{-1} (z)$ is known as the Jost function. The recurrent procedure applied to equations \e{eq:ab} shows that
\begin{equation}
2^{N+1}\omega (z) a_{0}\cdots a_{N-1} u_{-1}(z)= (1-4 a_{N-1}^2) \omega^{2N}  +\sum_{k=1}^{2N-1} c_{k}\omega^k +1
\label{eq:ab1}\end{equation}
where the coefficients $c_{k}$ can, in principle, be calculated in terms of the coefficients $a_{0},\cdots, a_{N-1}$, 
$b_{0},\cdots, b_{N-1}$.


   As is well-known, the Jost function is linked to  the perturbation determinant $D (z) $ for the pair $H_{1}$, $H $ by the relation
\begin{equation}
D (z):=\det \big(I+ V R_{1}(z)\big) =2^{N+1}\omega a_{0}\cdots a_{N-1} u_{-1}(z).
\label{eq:ab2}\end{equation}
Comparing \e{eq:ab1} and  \e{eq:ab2}, we get the following result.

 \begin{proposition}\label{dd}
 Let a Jacobi  operator $H$ be given by formula  \e{eq:J} where $a_{n}=1/2$   and  $b_{n}= 0$ for $n\geq N$.
Then the perturbation determinant $D(z)$ for the pair $H_{1}$, $H$  is a polynomial of $\omega =\omega(z)$ of degree at most $2N$ and the coefficient at $\omega^{2N} $ is necessary strictly smaller than $1$.
 \end{proposition}

  Note also  that the constant term of $D(z)$  equals $1$, but this fact is trivial because $D(z)\to 1$ as $|z|\to\infty$ or equivalently $\omega\to 0$.

 To construct the resolvent $R(z)=(H-z I)^{-1}$ of the operator $H$, one introduces the solution $\varphi_{n} (z)$ of equations \e{eq:ab} satisfying the boundary conditions $\varphi_{-1} (z)=0$, $\varphi_0 (z)=1$. Clearly, $\varphi_{n} (z)$ is a polynomial of $z$ of degree $n$. It is well-known that
\begin{equation}
(R (z)e_{n},e_{m}) = u_{-1}(z)^{-1}\varphi_{n} (z) u_m(z), \q n\leq m,
\label{eq:ab3}\end{equation}
and $(R (z)e_m, e_n)=(R (z)e_{n},e_{m})$. Proofs of formulas \e{eq:ab2} and \e{eq:ab3} can be found in the book
\cite{Teschl} where Jacobi operators were considered in the space $\ell^2 ({\Bbb Z})$; the case of the operators in 
$\ell^2 ({\Bbb Z}_{+})$ is quite similar.  On the other hand, these formulas are basically the same as the corresponding formulas for second order differential  operators on the half-axis with coefficients of compact support.

\bigskip

   {\bf 7.2.} 
   Let us now discuss a particular case $N=1$ when
     \[
V  = b_{0} (\cdot, e_{0}) e_{0}+ (a_{0} - 1/2) \big((\cdot, e_{0}) e_{1}  +  (\cdot, e_{1}) e_0 \big).
\]
In accordance with the notation of Section~1, we set $a=2a_{0}$, $b=2b_{0}$ and $V=V_{a,b}$. The operator $V_{a,b} $ has rank three so that the scheme of the main part of the paper can be easily adjusted to handle  the operator    $H_{a,b} =H_{1}+V_{a,b}$. Similarly, given that we have already obtained the explicit formulas for all matrix elements $(R_{a} (z) e_{n}, e_{m})$, we can consider  $H_{a,b}$ as a rank one perturbation of the operator $H_{a}$.

    The   generalized Chebyshev polynomials $\Ch_n(z;a,b)$
     are defined by the relations   $\Ch_0(z;a,b)=1$, $\Ch_1(z;a,b)=a^{-1} (2z-b)$ and \e{eq:chab}. They coincide of course with the functions $\varphi_{n} (z)$.   Formula \e{eq:CH2}  for $\Ch_n(z;a,b)$ remains true if function \e{eq:CH1} is replaced by a more general function
        \[
\gamma_{\pm}  (z;a,b)= \frac{a}{2} \pm  \frac{(a^2-2)z +b}{2a\sqrt{z^2-1}}.
\]


       To calculate the perturbation determinant for the pair $H_{1}, H_{a,b}$, we  
  here  use the approach described in the previous subsection.       The Jost solution is given by the equalities $u_{n} (z)= \omega(z)^n$ for $n\geq 1$, $u_{0} =a^{-1}$ and
   $2 u_{-1}(z)=(2z -b) u_{0} - a u_1 (z)$.
Therefore relation  \e{eq:ab2} yields the  formula
   \[
D_{a,b} (z)=(1- a ^2) \omega^{2}  -b \omega +1 .
\]
 Of course this formula reduces to
 \e{eq:D} if $b=0$. Thus we arrive at the following simple result.
 
 \begin{proposition}\label{d1}
  Let a Jacobi  operator $H$ be given by formula  \e{eq:J} where $a_{n}=1/2$   and  $b_{n}= 0$ for $n\geq 1$.
Then a polynomial 
\[
L(z)= l_{2}\omega^2+ l_{1}\omega+  l_0, \q \omega=\omega(z),
\]
is the perturbation determinant for the pair $H_{1}$, $H$   if and only if
$l_{2} < 1$, $l_{0}=1$ $($the coefficient $l_{1}$ remains arbitrary$)$. In this case the coefficients of $H$ can be recovered by the formulas $a_{0}=2^{-1}\sqrt{1-l_{2}}$ and $b_{0}=-l_1 /2$.
 \end{proposition}
 
 It is also   easy to find explicitly the perturbation determinant for $N=2$:
 \[
 D(z)= \alpha_{1}\omega^4 -2 (\alpha_{1} b_{0}+b_{1})\omega^3+ (\alpha_{1}+\alpha_{0}+4b_{0}b_{1})\omega^2
 -2 (b_{0}+b_{1})\omega+ 1, \q \omega=\omega(z),
 \]
 where for short we have put $\alpha_{j}=1-4 a_{j}^2$, $j=1,2$.
  Let us state   a generalization of Proposition~\ref{d1} to $N=2$.
 
  \begin{proposition}\label{d2}
        Let a Jacobi  operator $H$ be given by formula  \e{eq:J} where $a_{n}=1/2$   and  $b_{n}= 0$ for $n\geq 2$.
Then a polynomial 
\[
L(z)=\ell_4\omega^4 + \ell_3\omega^3+ \ell_{2}\omega^2+ \ell_{1}\omega+  \ell_0, \q \omega=\omega(z),
\]
is the perturbation determinant for the pair $H_{1}$, $H$  if and only if
$\ell_4 < 1$, $\ell_{0}=1$, the coefficients $\ell_{1}$ and $\ell_3$ are arbitrary and
\[
\ell_2<1+ \ell_4 +\frac{(\ell_3-\ell_1) (\ell_{1} \ell_4-\ell_{3}) }{(1-\ell_{4})^2}.
\]
  In this case the coefficients of $H$ can be recovered by the formulas  $a_1=2^{-1}\sqrt{1-\ell_4}$, 
\[
 b_{0}=  \frac{\ell_{3}-\ell_{1}}{2(1-\ell_4)},\q  b_{1}=  \frac{ \ell_{1} \ell_4 -\ell_{3}}{2(1-\ell_4)} 
\]
and
\[
     a_{0} =2^{-1}(1-\ell_4)^{-1}\sqrt{(1-\ell_2 +\ell_4) (1-\ell_4)^2+ ( \ell_{3}-\ell_{1}) (\ell_{1} \ell_4 -\ell_{3}) }.
 \]
  \end{proposition}
 
 Unfortunately, for $N\geq 3$, the results of this type become quite messy. Note that    the perturbation determinants were described for all $N$  in the paper \cite{DS} where however the results were stated in a less explicit form (not in terms of coefficients of polynomials).
 
 \appendix{}
 
 \section{Back to general Jacobi operators}

 Here we   return to  arbitrary bounded Jacobi operators \e{eq:J} and  discuss two general facts  which we were   unable to find in the literature.
 
 \bigskip

 {\bf A.1.}
 In the paper \cite{GS}, the authors posed the problem of characterizing the spectral measures $d\rho(\lambda)$ such that, for the corresponding Jacobi operator $H$, all coefficients $a_{n}=1$, that is, $H$ is a discrete Schr\"odinger operator. We will show that this problem admits a trivial solution in terms of   moments \e{eq:Ha2} of the measure $d\rho(\lambda)$. Let us set $h_{0}=1$ and
 \begin{equation*}
h_{n}= \det
\begin{pmatrix}
\kappa_{0}& \kappa_1 &  \cdots & \kappa_{n-1}\\
\kappa_{1}&  \kappa_2 &  \cdots & \kappa_{n}\\
      \vdots&    \vdots &  \vdots &  \vdots\\
\kappa_{n-1}&  \kappa_n &  \cdots & \kappa_{2n-2}
\end{pmatrix} 
 \end{equation*}
for $n\geq 1$.  Here the moments $\kappa_{0}=1, \kappa_1,\ldots, \kappa_{2n-2}$ are defined by \e{eq:Ha2}.

 \begin{proposition}\label{A2}
For a Jacobi operator \e{eq:J}, all coefficients $a_{n}=1$ if and only if $h_{n}=1$ for all $n\in{\Bbb Z}_{+}$.
 \end{proposition}
 
 \begin{pf}
 Let us proceed from the well-known formula (see, e.g., Theorem~A.2 in \cite{Si})
\[
a_{n}^2=  h_{n} h_{n+2}h^{-2}_{n+1}.
\]
Thus $a_{n}=1$ if $h_{n}=1$ for all $n$. Conversely, if $a_{n}=1$   for all $n$, then
\[
h_{n+2}=   h_{n+1}^2  h_{n}^{-1}.
\]
Since $h_{0}=h_{1}=1$, we find successively that $h_n=1$ for all $n$.
  \end{pf}

  \bigskip
  
 {\bf A.2.}
 Next, we consider the opposite situation and  characterize the spectral measures $d\rho(\lambda)$ of Jacobi operators $H$ with the coefficients $b_{n}=0$ for all $n\in{\Bbb Z}_{+}$. 
 
   \begin{proposition}\label{A1}
       Let a Jacobi  operator $H$ be given by formula  \e{eq:J}.
   Then coefficients $b_{n}=0$ for all $n\in{\Bbb Z}_{+}$   if and only if the spectral measure $d\rho(\lambda)$ of $H$ is even, i.e., $\rho(X)=\rho(-X)$ for all Borelian sets $X\subset{\Bbb R}$.
 \end{proposition}
 
 \begin{pf}
 Let the operators $J,B: \ell^2 ({\Bbb Z}_{+})\to \ell^2 ({\Bbb Z}_{+})$ be defined by the formulas $(Jf)_{n}= (-1)^n f_{n}$ and
 $(B f)_{n}= b_{n} f_{n}$. Then
   \[
(- JHJ f)_{n}= a_{n-1} f_{n-1}-b_n f_n+a_n f_{n+1}
\]
so that $-JHJ$ is a Jacobi operator and the condition $b_{n}=0$ for all $n\in{\Bbb Z}_{+}$ is equivalent to the relation
    \begin{equation}
 -JH J= H.
\label{eq:JJ2}\end{equation}

In terms of the Weyl function
\[
((H-zI)^{-1}e_{0}, e_{0})=\int_{-\infty}^\infty (\lambda-z)^{-1} d\rho(\lambda) ,
\]
the condition that the measure $d\rho(\lambda)$ is even, i.e., $d\rho(-\lambda)=- d\rho (\lambda)$, is equivalent to the relation
    \begin{equation}
((H+zI)^{-1}e_{0}, e_{0})=-((H-zI)^{-1}e_{0}, e_{0}).
\label{eq:JJ3}\end{equation}
Observe that
    \[
((H+zI)^{-1}e_{0}, e_{0})= ((H+zI)^{-1}Je_{0}, Je_{0})=((J H J+zI)^{-1} e_{0}, e_{0}),
\]
and hence \e{eq:JJ3} is equivalent to the relation
\[
((-J H J-zI)^{-1} e_{0}, e_{0}) =((H-zI)^{-1}e_{0}, e_{0}).
\]
Since the Weyl functions of Jacobi operators coincide if and only if   these operators coincide, the last relation is equivalent to  \e{eq:JJ2}.
   \end{pf}

\end{document}